\documentclass[5p,times,authoryear]{elsarticle}

\usepackage{setspace}

\usepackage{lineno,hyperref}
\modulolinenumbers[5]

\usepackage{lipsum} 
\usepackage{amsmath,amsthm,amsfonts, amssymb}
\usepackage{bbm} 
\usepackage{subcaption}
\usepackage{graphicx}
\graphicspath{{./figures/}}
\usepackage{epstopdf}

\usepackage{multicol}

\usepackage{algorithm,algorithmic}%

\usepackage{tikz,times,adjustbox}
\usetikzlibrary{mindmap,trees,backgrounds}
\usepackage{color}
\definecolor{red1}{rgb}{0.502,0,0}
\definecolor{blue1}{rgb}{0.098,0.098,0.498}
\definecolor{brown1}{rgb}{0.545,0.271,0.075}
\definecolor{green1}{rgb}{0,0.392,0}

\definecolor{red2}{rgb}{0.698,0.1333,0.1333}
\definecolor{blue2}{rgb}{0.275,0.51,0.706} 
\definecolor{brown2}{rgb}{0.824,0.412,0.118} 
\definecolor{green2}{rgb}{0.180,0.545,0.341}

\definecolor{red3}{rgb}{0.8039,0.3608,0.3608}
\definecolor{blue3}{rgb}{0.1176,0.565,1} 
\definecolor{brown3}{rgb}{1.0000,0.5490,0} 
\definecolor{green3}{rgb}{0.2353,0.7020,0.4431}

\theoremstyle{plain}
\newtheorem{thm}{Theorem}
\newtheorem{lem}{Lemma}
\newtheorem{prop}{Proposition}
\newtheorem{cor}{Corollary}
\newtheorem{assumption}{Assumption}
\theoremstyle{definition}
\newtheorem{defn}{Definition}

\theoremstyle{remark}
\newtheorem{rem}{Remark}

\DeclareMathOperator*{\esssup}{ess\,sup}
\DeclareMathOperator{\VaR}{VaR}
\DeclareMathOperator{\CVaR}{CVaR}
\DeclareMathOperator{\diag}{diag}
\newcommand{\floor}[1]{\left\lfloor #1 \right\rfloor}
\newcommand{\ceil}[1]{\left\lceil #1 \right\rceil}

\usepackage{multirow,array,float}

\newcolumntype{P}[1]{>{\centering\arraybackslash}p{#1}}
\newcolumntype{L}[1]{>{\flushleft\arraybackslash}p{#1}}

\begin{document}
\sloppy
{\onecolumn \large 
\section*{About} 
\begin{itemize}
	\item This is Part 1 of a two-part review paper titled ``Data-driven Decision Making in Power Systems with Probabilistic Guarantees: Theory and Applications of Chance-constrained Optimization'' by Xinbo Geng and Le Xie, Annual Reviews in Control (under review).
	\item Part 1 ``Data-driven Decision Making with Probabilistic Guarantees (Part I): A Schematic Overview of Chance-constrained Optimization'' is available at arXiv:1903.10621.
	\item Part 2 ``Data-driven Decision Making in with Probabilistic Guarantees (Part II): Applications of Chance-constrained Optimization in Power Systems'' is available at arXiv.
	\item The Matlab Toolbox \emph{ConvertChanceConstraint} (CCC) is available at \url{https://github.com/xb00dx/ConvertChanceConstraint-ccc}.
\end{itemize}
Please let us know if we missed any critical references or you found any mistakes in the manuscript.
\section*{Recent Updates} 
\begin{description}
	\item [04/2019] More CVaR-based (Convex Approximation) results are added in Part 1.
	\item [02/2019] Toolbox published at \url{https://github.com/xb00dx/ConvertChanceConstraint-ccc}. We are still working on the toolbox website and documents.
\end{description}
}

\begin{frontmatter}

\title{Data-driven Decision Making with Probabilistic Guarantees (Part I):\\
A Schematic Overview of Chance-constrained Optimization}

\author{Xinbo Geng, Le Xie}
\address{Texas A\&M University, College Station, TX, USA.}




\begin{abstract}
Uncertainties from deepening penetration of renewable energy resources have posed critical challenges to the secure and reliable operations of future electric grids.
Among various approaches for decision making in uncertain environments, this paper focuses on chance-constrained optimization, which provides explicit probabilistic guarantees on the feasibility of optimal solutions.
Although quite a few methods have been proposed to solve chance-constrained optimization problems, there is a lack of comprehensive review and comparative analysis of the proposed methods. Part I of this two-part paper reviews three categories of existing methods to chance-constrained optimization: (1) scenario approach; (2) sample average approximation; and (3) robust optimization based methods. Data-driven methods, which are not constrained by any particular distributions of the underlying uncertainties, are of particular interest. Part II of this two-part paper provides a literature review on the applications of chance-constrained optimization in power systems. Part II also provides a critical comparison of existing methods based on numerical simulations, which are conducted on standard power system test cases.
\end{abstract}

\begin{keyword}
data-driven, power system, chance constraint, probabilistic constraint, stochastic programming, robust optimization, chance-constrained optimization.
\end{keyword}

\end{frontmatter}


\section{Introduction}

The objective of this article is to provide a comprehensive and up-to-date review of mathematical formulations, computational algorithms, and engineering implications of chance-constrained optimization in the context of electric power systems. 
In particular, this paper focuses on the data-driven approaches to solving chance-constrained optimization without knowing or making specific assumptions on the underlying distribution of the uncertainties. A more general class of problem, i.e. distributionally robust optimization or ambiguous chance constraint, is beyond the scope of this paper.

\subsection{An Overview of Chance-constrained Optimization} 
\label{sub:an_overview_of_chance_constrained_optimization}
Chance-constrained optimization (CCO) is an important tool for decision making in uncertain environments. Since its birth in 1950s, CCO has found many successful applications in various fields, e.g. economics \citep{yaari_uncertain_1965}, control theory \citep{calafiore_scenario_2006}, chemical process \citep{sahinidis_optimization_2004,henrion_stochastic_2001}, water management \citep{dupacova_stochastic_1991} and recently in machine learning \citep{xu_robust_2009,ben-tal_robust_2009,caramanis_14_2012,ben-tal_chance_2011,sra_optimization_2012,gabrel_recent_2014}. Chance-constrained optimization plays a particularly important role in the context of electric power systems \citep{ozturk_solution_2004,wang_chance-constrained_2012}, applications of CCO can be found in various time-scales of power system operations and at different levels of the system. 

The first chance-constrained program was formulated in \citep{charnes_cost_1958},
then was extensively studied in the following 50 years, e.g. \citep{charnes_chance-constrained_1959,charnes_deterministic_1963,kataoka_stochastic_1963,pinter_deterministic_1989,sen_relaxations_1992,prekopa_programming_1998,ruszczynski_stochastic_2003-1,ben-tal_robust_2009,prekopa_stochastic_1995}. Previously, most methods to solve CCO problems deal with specific families of distributions, such as log-concave distributions \citep{miller_chance_1965,prekopa_stochastic_1995}. Many novel methods appeared in the past ten years, e.g. scenario approach \citep{calafiore_scenario_2006}, sample average approximation \citep{ruszczynski_probabilistic_2002,luedtke_sample_2008} and convex approximation \citep{nemirovski_convex_2006}. Most of them are generic methods that are not limited to specific distribution families and require very limited knowledge about the uncertainties. In spite of many successful applications of these methods in various fields, there is a lack of comprehensive review and a critical comparison.

\subsection{Contributions} 
\label{sub:contributions_of_this_paper}
The main contributions of this paper are twofold:
\begin{enumerate}
	\item We provide a detailed tutorial on existing algorithms to solve chance-constrained programs and a survey of major theoretical results. To the best of our knowledge, there is no such review available in the literature;
	\item We implement all the reviewed methods and develop an open-source Matlab toolbox (ConvertChanceConstraint), which is available on Github \footnote{github.com/xb00dx/ConvertChanceConstraint-ccc}.
\end{enumerate}

\subsection{Organization of This Paper} 
\label{sub:organization_of_this_paper}
The remainder of this paper is organized as follows. Section \ref{sec:chance_constrained_optimization} introduces chance-constrained optimization. Section \ref{sec:properties_chance_constrained_optimization} summarizes the fundamental properties of chance-constrained optimization problems. An overview of how to solve the chance-constrained optimization problem is described in Section \ref{sec:an_overview_of_solutions_to_cco}, which outlines Section \ref{sec:scenario_approach}-\ref{sec:robust_optimization}. Three major approaches to solving chance-constrained optimization (scenario approach, sample average approximation and robust optimization based methods) are presented in Section \ref{sec:scenario_approach}-\ref{sec:robust_optimization}, respectively. The structure and usage of the Toolbox \emph{ConvertChanceConstraint} is in Section \ref{sec:convertchanceconstraint_a_matlab_toolbox}. Concluding remarks are in Section \ref{sec:concluding_remarks}.

\subsection{Notations} 
\label{sub:notations}
The notations in this paper are standard. All vectors and matrices are in the real field $\mathbf{R}$. Sets are in calligraphy fonts, e.g. $\mathcal{S}$. The upper and lower bounds of a variable $x$ are denoted by $\overline{x}$ and $\underline{x}$. The estimation of a random variable $\epsilon$ is $\hat{\epsilon}$. We use $\mathbf{1}_n$ to denote an all-one vector in $\mathbf{R}^n$, the subscript $n$ is sometimes omitted for simplicity. The absolute value of vector $x$ is $|x|$, and the cardinality of a set $\mathcal{S}$ is $|\mathcal{S}|$. Function $[a]_+$ returns the positive part of variable $a$. The indicator function $\mathbbm{1}_{x>0}$ is one if $x > 0$. The floor function $\floor{a}$ returns the largest integer less than or equal to the real number $a$. 
The ceiling function $\ceil{a}$ returns the smallest integer greater than or equal to $a$.
$\mathbb{E}[\xi]$ is the expectation of a random vector $\xi$, $\mathbb{V}(x)$ denotes the violation probability of a candidate solution $x$, and $\mathbb{P}_\xi(\cdot)$ is the probability taken with respect to $\xi$.
The transpose of a vector $a$ is $a^\intercal$. Infimum, supremum and essential supremum are denoted by $\inf$, $\sup$ and $\esssup$. The element-wise multiplication of the same-size vectors $a$ and $b$ is denoted by $a \circ b$.

\section{Chance-constrained Optimization} 
\label{sec:chance_constrained_optimization}
\subsection{Introduction} 
\label{sub:introduction}
We study the following chance-constrained optimization problem throughout this paper:
\begin{subequations}
\label{form:cco_joint}
\begin{align}
  \text{(CCO):}~\min_{x}~&  c^\intercal x  \\
  \text{s.t.}~&  \mathbb{P}_{\xi} \Big( f(x, \xi) \le 0 \Big) \ge 1 - \epsilon \label{form:ccp_joint_cc} \\
  & x \in \mathcal{X} 
\end{align}
\end{subequations}
where $x \in \mathbf{R}^n$ is the decision variable and random vector $\xi \in \mathbf{R}^d$ is the source of uncertainties. Without loss of generality \footnote{Using the epigraph formulation as mentioned in \citep{campi_scenario_2009,boyd_convex_2004}.}, we assume the objective function is linear in $x$ and does not depend on $\xi$.
Constraint (\ref{form:ccp_joint_cc}) is the \emph{chance constraint} (or \emph{probabilistic constraint}), it requires the inner constraint $f(x,\xi) \le 0$ to be satisfied with high probability $1- \epsilon$. The inner constraint $f(x,\xi): \mathbf{R}^n \times \mathbf{R}^d \rightarrow \mathbf{R}^m$ consists of $m$ individual constraints, i.e. $f(x, \xi) = \big(f_1(x,\xi),f_2(x,\xi),\cdots,f_m(x,\xi) \big)$. Set $\mathcal{X}$ represent the deterministic constraints. Parameter $\epsilon$ is called the \emph{violation probability} of (CCO). Notice that $f(x,\xi)$ is random due to $\xi$, the probability $\mathbb{P}$ is taken with respect to $\xi$. Sometimes the probability is denoted by $\mathbb{P}_\xi$ to avoid confusion.

It is worth mentioning that CCO is closely related with the theory of risk management. For example, an individual chance constraint $\mathbb{P}(f_i(x,\xi) \le 0) \ge 1- \epsilon_i$ can be equivalently interpreted as a constraint on the value at risk $\VaR(f_i(x,\xi); 1- \epsilon_i) \le 0$. This connection can be directly seen from the definition.
\begin{defn}[Value at Risk]
Value at risk (VaR) of random variable $\zeta$ at level $1- \epsilon$ is defined as 
\begin{equation}
	\VaR(\zeta;1- \epsilon) := \inf\big\{\gamma: \mathbb{P}(\zeta \le \gamma ) \ge 1- \epsilon \big\}
\end{equation}
\end{defn}
More details about this can be found in Section \ref{ssub:convex_approximation}, \citep{rockafellar_optimization_2000,chen_cvar_2010} and references therein.

CCO is closely related with two other major tools for decision making with uncertainties: stochastic programming and robust optimization. The idea of sample average approximation, which originated from stochastic programming, can be applied on chance-constrained programs (Section \ref{sec:sample_average_approximation}). Section \ref{sec:robust_optimization} demonstrates the connection between robust optimization and CCO.

\subsection{Joint and Individual Chance Constraints} 
\label{sub:joint_and_individual_chance_constraints}
Constraint (\ref{form:ccp_joint_cc}) is called a \emph{joint chance constraint} because of its multiple inner constraints \citep{miller_chance_1965}, i.e. 
\begin{equation}
\label{eqn:joint_chance_constraint}
\mathbb{P}\Big(f_1(x, \xi) \le 0, f_2(x, \xi) \le 0, \cdots,f_m(x, \xi) \le 0\Big) \ge 1 -\epsilon
\end{equation}
Alternatively, each one of the following $m$ constraints is called an \emph{individual chance constraint}:
\begin{equation}
\label{eqn:individual_chance_constraints}
\mathbb{P}\Big( f_i(x,\xi) \le 0\Big) \le 1 - \epsilon_i,~i = 1,2,\cdots,m
\end{equation}
Joint chance constraints typically have more modeling power since an individual chance constraint is a special case ($m=1$) of a joint chance constraint. But individual chance constraints are relatively easier to deal with (see Section \ref{sub:special_cases} and \ref{sub:safe_approximation_of_individual_chance_constraints}). There are several ways to convert individual and joint chance constraints between each other. 

First, a joint chance constraint can be written as a set of individual chance constraints using Bonferroni inequality or Boole's inequality. 
Notice (\ref{eqn:joint_chance_constraint}) can be represented as
\begin{equation}
	\mathbb{P}_\xi\Big(\cup_{i=1}^m \big\{f_i(x,\xi) \ge 0\big\}\Big) \le \epsilon.
\end{equation}
Since $\mathbb{P}_\xi(\cup_{i=1}^m \{f_i(x,\xi) \ge 0\}) \le \sum_{i=1}^m \mathbb{P}_\xi(\{f_i(x,\xi) \ge 0\})$, if $\sum_{i=1}^{m} \epsilon_i \le \epsilon$,  then any feasible solution to (\ref{eqn:individual_chance_constraints}) is also feasible to (\ref{eqn:joint_chance_constraint}). In other words, (\ref{eqn:individual_chance_constraints}) is a \emph{safe approximation} (see Definition \ref{defn:safe_approx}) to (\ref{eqn:joint_chance_constraint}) when $\sum_{i=1}^{m} \epsilon_i \le \epsilon$.
With appropriate $\{\epsilon_i\}_{i=1}^m$, (\ref{eqn:individual_chance_constraints}) could be a good approximation of (\ref{eqn:joint_chance_constraint}). However, it is usually difficult to find such $\{\epsilon_i\}_{i=1}^m$. Some other issues of this approach are discussed in Section \ref{ssub:converting_a_joint_chance_constraint_to_individual}.

Alternatively, a joint chance constraint (\ref{eqn:joint_chance_constraint}) is equivalent to the following individual chance constraint:
\begin{equation}
\label{eqn:individual_cc_equivalent_joint}
\mathbb{P}_\xi\big( \overline{f}(x,\xi) \le 0 \big) \ge 1 -\epsilon
\end{equation}
where $\overline{f}(x,\xi):\mathbf{R}^{n} \times \mathbf{R}^d \rightarrow \mathbf{R}$ is the pointwise maximum of functions $\{f_i(x,\xi)\}_{i=1}^m$ over $x$ and $\xi$, i.e. 
\begin{equation}
\overline{f}(x,\xi):=\max\Big\{f_1(x,\xi),f_2(x,\xi),\cdots,f_m(x,\xi) \Big\}.
\end{equation}
It is worth noting that converting $\{f_i(x,\xi)\}_{i=1}^m$ to $\overline{f}(x,\xi)$ could lose nice structures of the original constraint $f(x,\xi) \le 0$ and cause more difficulties. 

In this paper, we focus on the chance-constrained optimization problems with a \emph{joint} chance constraint.

\subsection{Critical Definitions and Assumptions} 
\label{sub:critical_definitions_assumptions}
Theoretical results in the following sections are based on the critical definitions and assumptions below.
\begin{defn}[Violation Probability]
Let $x^\diamond$ denote a candidate solution to (CCO), its violation probability is defined as 
\begin{equation}
  \mathbb{V}(x^\diamond) := \mathbb{P}_{\xi} \Big( f(x^\diamond, \xi) \ge 0 \Big)
\end{equation}
\end{defn}
\begin{defn}
\label{defn:cco-feasible}
$x^\diamond$ is a \emph{feasible} solution for (CCO) if $x^\diamond \in \mathcal{X}$ and $\mathbb{V}(x^\diamond) \le \epsilon$. Let $\mathcal{F}_\epsilon$ denote the set of feasible solutions to the chance constraint (\ref{form:ccp_joint_cc}),
\begin{equation*}
	\mathcal{F}_\epsilon:= \{x \in \mathbf{R}^n: \mathbb{V}(x) \le \epsilon \} = \{x \in \mathbf{R}^n: \mathbb{P}_\xi\Big(f(x,\xi) \le 0\Big) \ge 1- \epsilon \}
\end{equation*}
then $x^\diamond$ is \emph{feasible} to (CCO) if $x^\diamond \in \mathcal{X} \cap \mathcal{F}_\epsilon$.
\end{defn}
Although (CCO) seeks optimal solutions under uncertainties, it is a \emph{deterministic} optimization problem. To better see this, (CCO) can be equivalently written as $\min_{x \in \mathcal{X}}~ c^\intercal x,~\text{s.t.}~\mathbb{V}(x) \le \epsilon$ or $\min_{x \in \mathcal{X} \cap \mathcal{F}} c^\intercal x$.
\begin{defn}
Let $o^\star$ denote the optimal objective value of (CCO). For simplicity, we define $o^\star = +\infty$ when (CCO) is infeasible and $o^\star = -\infty$ when (CCO) is unbounded. Let $x^\star$ denote the optimal solution to (CCO) if exists, and $o^\star = c^\intercal x^\star$. 
\end{defn}
\begin{defn}
We say a candidate solution $x^\diamond$ is \emph{conservative} if $\mathbb{V}(x^\diamond) \ll \epsilon$ or $ c^\intercal x^\diamond \gg o^\star$.
\end{defn}
Most existing theoretical results on (CCO) are built upon the following two assumptions.
\begin{assumption}
\label{assu:distribution}
Let $\Xi$ denote the support of random variable $\xi$, the distribution $\xi \sim \Xi$ exists and is fixed.
\end{assumption}
Assumption \ref{assu:distribution} only assumes the existence of an underlying distribution, but we do not necessarily need to know it to solve (CCO). Removing assumption \ref{assu:distribution} leads to a more general class of problem named \emph{distributionally robust optimization} or \emph{ambiguous chance constraints}.
Section \ref{sub:ambiguous_chance_constraints} discusses cases with Assumption \ref{assu:distribution} removed.
\begin{assumption}
\label{assu:convexity}
(1) Function $f(x, \xi)$ is convex in $x$ for every instance of $\xi$, and (2) the deterministic constraints define a convex set $\mathcal{X}$.
\end{assumption}
The convexity assumption makes it possible to develop theories on (CCO). However, (CCO) is often non-convex even under Assumption \ref{assu:convexity}. More details are presented in Section \ref{sub:hardness} and difficulty (D2).



\section{Fundamental Properties} 
\label{sec:properties_chance_constrained_optimization}
\subsection{Hardness} 
\label{sub:hardness}
Although CCO is an important and useful tool for decision making under uncertainties, it is very difficult to solve in general. Major difficulties come from two aspects: 
\begin{description}
	\item [(D1)] It is difficult to check the feasibility of a candidate solution $x^\diamond$. Namely, it is intractable to evaluate the probability $\mathbb{P}_\xi\big(f(x^\diamond,\xi) \le 0\big)$ with high accuracy. More specifically, calculating probability involves multivariate integration, which is NP-Hard \citep{khachiyan_problem_1989}. The only general method might be Monte-Carlo simulation, but it can be computationally intractable due to the curse of dimensionality.
	\item [(D2)] It is difficult to find the optimal solution $x^\star$ and $o^\star$ to (CCO). Even with the convexity assumption (Assumption \ref{assu:convexity}), the feasible region of (CCO) is often non-convex except a few special cases. For example, Section \ref{sub:feasible_region} shows the feasible region of (CCO) with separable chance constraints is a union of cones, which is non-convex in general. Although researchers have proved various sufficient conditions on the convexity of (CCO), it remains challenging to solve (CCO) because of difficulty (D1). Most of times, however, we are agnostic about the properties of the feasible region $\mathcal{F}_\epsilon$.
\end{description} 
Despite that fact that Assumptions \ref{assu:distribution} and \ref{assu:convexity} greatly simplify the problem and make theoretical analysis on (CCO) possible, (D1) and (D2) still exist and pose great challenges to solve (CCO).
\begin{thm}[\citep{luedtke_integer_2010,qiu_covering_2014}]
\label{thm:cco_np_hardness}
(CCO) is strongly NP-Hard.
\end{thm}
\begin{thm}[\citep{ahmed_relaxations_2018-1}]
\label{thm:cco_np_hard_to_approx}
Unless $\text{P}=\text{NP}$, it is impossible to obtain a polynomial time algorithm for (CCO) with a constant approximation ratio.
\end{thm}
Theorem \ref{thm:cco_np_hardness} formalizes the hardness results of solving (CCO), Theorem \ref{thm:cco_np_hard_to_approx} further demonstrates it is also difficult to obtain approximate solutions to (CCO): any polynomial algorithm is not able to find a solution $x^\ast$ (with $o^\ast = c^\intercal x^\ast$) such that $|o^\ast/o^\star|$ is bounded by a constant $C$. In other words, any polynomial-time algorithm could be arbitrarily worse. 

\subsection{Special Cases} 
\label{sub:special_cases}
Although (CCO) is NP-Hard to solve in general, there are several special cases in which solving (CCO) is relatively easy.
The most well-known special case is (\ref{form:cco_gaussian}), which was first proved in \citep{kataoka_stochastic_1963}.
\begin{subequations}
\label{form:cco_gaussian}
\begin{align}
\min_{x \in \mathcal{X}} ~ & c^\intercal x \\
\text{s.t.}~& \mathbb{P}( a^\intercal x + b^\intercal  \xi + \xi^\intercal D x \le e)  \ge 1- \epsilon \label{form:cco_gaussian_cc}
\end{align}
\end{subequations}
Parameters $a \in \mathbf{R}^n$,$b \in \mathbf{R}^d$, $D \in \mathbf{R}^{d \times n}$ and $e \in \mathbf{R}$ are fixed coefficients. 
$\xi \sim \mathcal{N}(\mu, \Sigma)$ is a multivariate Gaussian random vector with mean $\mu$ and covariance $\Sigma$. Notice that (\ref{form:cco_gaussian_cc}) is an \emph{individual} chance constraint with multivariate Gaussian coefficients. Let $\Phi(\cdot)^{-1}$ denote the inverse cumulative distribution function (CDF) function of a standard normal distribution.
It is easy to show that if $\epsilon \le 1/2$, (\ref{form:cco_gaussian}) is equivalent to (\ref{form:cco_gaussian_socp}), which is a second order cone program (SOCP) and can be solved efficiently.
\begin{subequations}
\label{form:cco_gaussian_socp}
\begin{align}
\min_{x \in \mathcal{X}} ~ & c^\intercal x \\
\text{s.t.}~& e - b^\intercal \mu - (a + D^\intercal \mu)^\intercal x  \ge \nonumber \\
& \qquad \Phi^{-1}(1- \epsilon ) \sqrt{ (b+Dx)^\intercal \Sigma (b+Dx)} \label{form:cco_gaussian_soc}	
\end{align}
\end{subequations}
The case of log-concave distribution \citep{prekopa_logarithmic_1971,prekopa_stochastic_1995,prekopa_uniform_2011} is another famous special case where chance constraint is convex. 
There are many other sufficient conditions on the convexity of chance constraints, e.g. \citep{lagoa_convexity_1999,calafiore_distributionally_2006,henrion_convexity_2008,henrion_convexity_2011,van_ackooij_eventual_2015}.

\subsection{Feasible Region} 
\label{sub:feasible_region}
A chance-constrained program with only right hand side uncertainties (\ref{form:cco_rhs_uncertain}) is considered in this section.
With this example, we provide deeper understandings on the non-convexity of (CCO).
\begin{subequations}
\label{form:cco_rhs_uncertain}
\begin{align}
\min_{x \in \mathcal{X}}~& c^\intercal x \\
\text{s.t.}~ & \mathbb{P}(f(x) \le \zeta) \ge p \label{form:cco_separable}
\end{align}
\end{subequations}
In (\ref{form:cco_separable}), the inner function $f(x): \mathbf{R}^n \rightarrow \mathbf{R}^m$ is deterministic. The only uncertainty  is the right-hand side value, represented by a random vector $\zeta \in \mathbf{R}^m$. Chance constraints like (\ref{form:cco_separable}) are also named \emph{separable} chance constraints (or probabilistic constraints) since the deterministic and random parts are separated. We replace $1- \epsilon$ with $p$ in (\ref{form:cco_separable}) to follow the convention in the existing literature.
\begin{defn}[$p$-efficient points \citep{shapiro_lectures_2009}]
Let $p \in (0,1)$, a point $v \in \mathbf{R}^m$ is called a \emph{$p$-efficient point} of the probability function $\mathbb{P}_\zeta(\zeta \le z)$, if $\mathbb{P}_\zeta(\zeta \le v) \ge p$ and there is no $z \le v$, and $z \ne v$ such that $\mathbb{P}_\zeta(\zeta \le z) \ge p$.
\end{defn}
\begin{thm}[\citep{shapiro_lectures_2009} \citep{prekopa_stochastic_1995}]
\label{thm:feasible_set_union_of_cones}
Let $\mathcal{E}$ be the index set of $p$-efficient points $v^i$, $i \in \mathcal{E}$. Let $\mathcal{F}_p := \{x \in \mathbf{R}^n: \mathbb{P}_\zeta(f(x) \le \zeta) \ge p\}$ denote the feasible region of (\ref{form:cco_separable}), then it holds that
\begin{equation}
	\mathcal{F}_p = \cup_{i \in \mathcal{E}} K_i
\end{equation}
where each cone $K_i$ is defined as $K_i := v^i + \mathbf{R}_+^m$, $i \in \mathcal{E}$.
\end{thm}
Theorem \ref{thm:feasible_set_union_of_cones} shows the geometric properties of (CCO). The finite union of convex sets need not to be convex, therefore the feasible region of (CCO) is generally non-convex.
\begin{rem}
Many methods to solve (CCO) (e.g. \citep{beraldi_probabilistic_2002,prekopa_programming_1998,kress_minmax_2007} ) start with a partial or complete enumeration of $p$-efficient points. However, the number of $p$-efficient points could be astronomic or even infinite. See \citep{shapiro_lectures_2009,prekopa_stochastic_1995} and references therein for the finiteness results of $p$-efficient points and complete theories and algorithms on $p$-efficient points.
\end{rem}

\subsection{Ambiguous Chance Constraints} 
\label{sub:ambiguous_chance_constraints}
\emph{Ambiguous chance constraint} is a generalization of chance constraints,
\begin{equation}
	\mathbb{P}_{\xi \sim P}\Big( f(x,\xi) \le 0 \Big) \ge 1 - \epsilon, ~\forall P \in \mathcal{P}.
\end{equation}
It requires the inner chance constraint $f(x,\xi) \le 0$ holds with probability $1- \epsilon$ for any distribution $P$ belonging to a set of pre-defined distributions $\mathcal{P}$.

Ambiguous chance constraints are particularly useful in the cases where only partial knowledge on the distribution $P$ is available, e.g. we know only that $P$ belongs a given family of $\mathcal{P}$. 
However, it is generally more difficult to solve ambiguous chance constraints, and the theoretical results rely on different assumptions of uncertainties. This paper only reviews solutions to CCO, studies on ambiguous chance constraints are beyond the scope of this paper.

\section{An Overview of Solutions to CCO} 
\label{sec:an_overview_of_solutions_to_cco}
This paper concentrates on solutions to (CCO) with the following properties: (i) dealing with both difficulties (D1) and (D2) mentioned in Section \ref{sub:hardness}; (ii) utilizing information from data (only) without making suspicious assumptions on the distribution of uncertainties; and (iii) possessing rigorous guarantees on the feasibility and optimality of returned solutions. Section \ref{sub:classification_of_solutions}-\ref{sub:theoretical_guarantees} explain these three properties in detail. Section \ref{sub:an_overview} provides an overview of methods with the properties above.

\subsection{Classification of Solutions} 
\label{sub:classification_of_solutions}
Existing methods on (CCO) can be roughly classified into four categories \citep{ahmed_solving_2008}:
\begin{description}
	\item[(C1)] When both difficulties (D1) and (D2) in Section \ref{sub:hardness} are absent, (CCO) is convex and the probability $\mathbb{P}(f(x,\xi)\le0)$ is easy to calculate.
	The only known case in this category is the individual chance constraint (\ref{form:cco_gaussian}) with Gaussian distributions, which might be the only special case of (CCO) that can be easily solved;
	\item[(C2)] When (D1) is absent but (D2) is present, it is relatively easy to calculate $\mathbb{P}(f(x,\xi)\le 0)$ (e.g. finite distributions with not too many realizations). As shown in Theorem \ref{thm:feasible_set_union_of_cones}, the feasible region of (CCO) could be non-convex and solutions typically rely on integer programming and global optimization \citep{ahmed_solving_2008};
	\item[(C3)] When (D1) is present but (D2) is absent, (CCO) is proved to be convex but remains difficult to solve because of the difficulty (D1) in calculating probabilities. This case often requires approximating the probability via simulations or specific assumptions. All examples mentioned in Section \ref{sub:special_cases} except (\ref{form:cco_gaussian}) belong to this category.
	\item[(C4)] When both difficulties (D1) and (D2) are present, it is almost impossible to find the optimal solution $x^\star$ and $o^\star$. All existing methods attempt to obtain approximate solutions or suboptimal solutions and construct upper and lower bounds on the true objective value $o^\star$ of (CCO).
\end{description}
Methods associated with (C1)-(C3) are briefly mentioned in Section \ref{sec:properties_chance_constrained_optimization}, the remaining part of this paper presents more general and powerful methods in category (C4).

\subsection{Prior Knowledge} 
\label{sub:prior_knowledge}
\begin{figure}[tb]
	\centering
	\includegraphics[width=\linewidth]{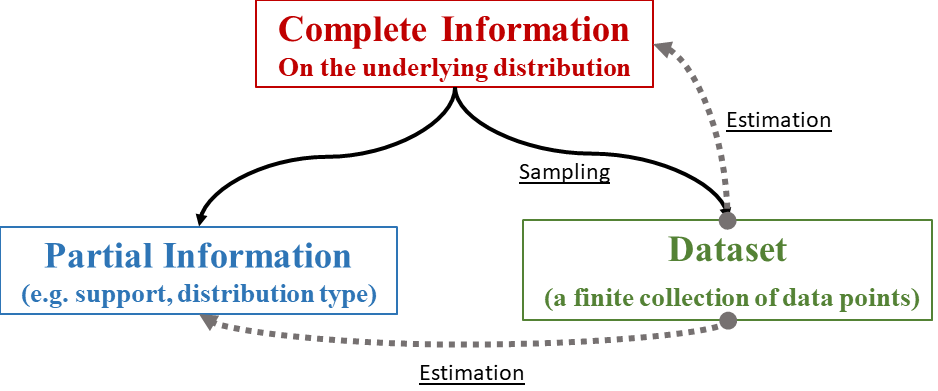}
	\caption{Different Knowledge Levels to Solve (CCO)}
	\label{fig:cco_information}
\end{figure}

In order to solve (CCO), a reasonable amount of prior knowledge on the underlying distribution $\xi \sim \Xi$ is necessary. Figure \ref{fig:cco_information} illustrates three categories of prior knowledge:
\begin{description}
	\item[(K1)] We know the exact distribution $\xi \sim \Xi$ thus have \emph{complete knowledge} on the underlying distribution;
	\item[(K2)] We know partially on the distribution (e.g. multivariate Gaussian distribution with bounded mean and variance) and thus have \emph{partial knowledge};
	\item[(K3)] We have a finite dataset $\{\xi^i\}_{i=1}^N$, this is another case of \emph{partial knowledge}. 
\end{description}
It can be seen that prior information in (K2) is a strict subset of (K1), also by sampling we can construct a dataset in (K3) from the exact distribution in (K1). It seems (K1) is the best starting point to solve (CCO). However, \emph{probability distributions are not known in practice, they are just models of reality and exist only in our imagination.} What exists in reality is \emph{data}. Therefore (K3) is the most practical case and becomes the focus of this paper. Almost all the data-driven methods to solve (CCO) are based on the following assumption.
\begin{assumption}
\label{assu:iid_samples}
The samples (scenarios) $\xi^i$ ($i=1,2,\cdots,N$) in the \emph{dataset} $\{\xi^i\}_{i=1}^N$ are independent and identically distributed (i.i.d.).
\end{assumption}

\subsection{Theoretical Guarantees} 
\label{sub:theoretical_guarantees}
In this paper, we concentrate on the theoretical aspects of the reviewed methods. In particular, we pay special attention to \emph{feasibility guarantees} and \emph{optimality guarantees}.

Given a candidate solution $x^\diamond$ to (CCO), the first and possibly most important thing is to check its \emph{feasibility}, i.e. if $\mathbb{V}(x^\diamond) \le \epsilon$.
Although (D1) demonstrates the difficulty in calculating $\mathbb{V}(x^\diamond)$ with high accuracy, there are various feasibility guarantees that either estimate $\mathbb{V}(x^\diamond)$ or provide upper bound on $\mathbb{V}(x^\diamond)$. The feasibility results can be classified into two categories: a-priori and a-posteriori guarantees. The \emph{a-priori} ones typically provide prior conditions on (CCO) and the dataset $\{\xi^i\}_{i=1}^N$, the feasibility of the corresponding solution $x^\diamond$ is guaranteed \emph{before} obtaining $x^\diamond$. Examples of this type include Corollary \ref{cor:prior_sample_complexity_fully_supported}, Theorem \ref{thm:prior_guarantee_helly_dim},\ref{thm:rlo_safe_approx} and \ref{thm:convex_approx_safe_approx}. As the name suggests, the \emph{a-posteriori} guarantees make effects \emph{after} obtaining $x^\diamond$. The a-posteriori guarantees are constructed based on the observations of the structural features associated with $x^\diamond$. Examples include Theorem \ref{thm:posterior_guarantee} and Proposition \ref{prop:check_feasibility_scenario_approach}.

Given a candidate solution $x^\diamond$ and the associated objective value $o^\diamond = c^\intercal x^\diamond$, another important question to be answered is about the \emph{optimality} gap $|o^\diamond - o^\star|$. Although finding $o^\star$ is often an impossible mission because of difficulty (D2), bounding from below on $o^\star$ is relatively easier. Sections \ref{sub:constructing_lower_bounds_via_scenario_approach} and \ref{sub:optimality_guarantees_of_sample_average_approximation} dedicate to algorithms of constructing lower bounds $\underline{o} \le o^\star$. 



\subsection{A Schematic Overview} 
\label{sub:an_overview}
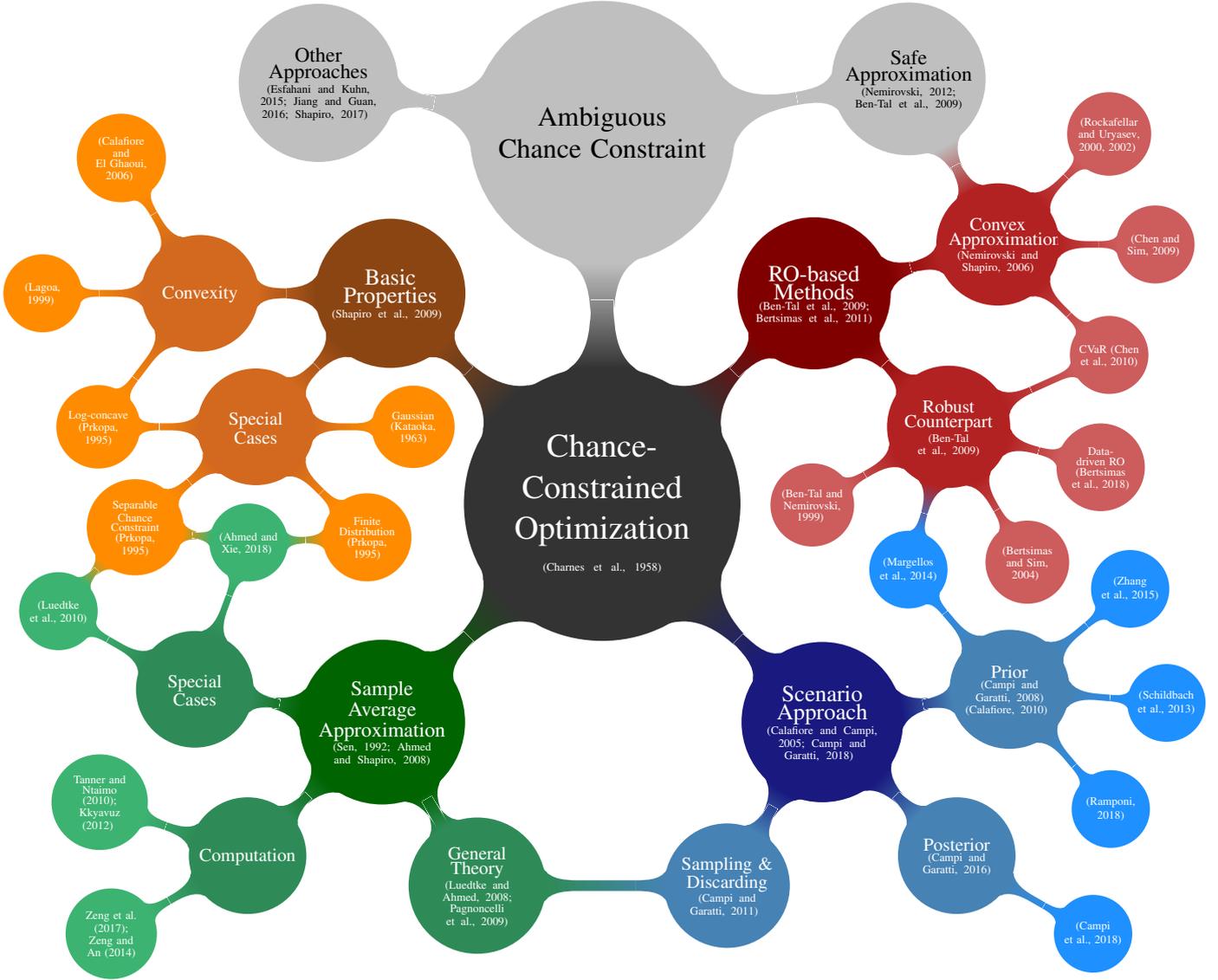
\begin{figure*}[tb]
	\centering
\begin{adjustbox}{width=\textwidth}	











\begin{tikzpicture}
  \begin{scope}[mindmap, concept color=lightgray,text=black]
    \node [concept] (ambiguous-cc) at (0,5.7) {Ambiguous Chance Constraint} 
      child [grow=10, level distance=135] 
        {node [concept] (safe-approx) { \begin{spacing}{0} Safe\\Approximation\\{\tiny\citep{nemirovski_safe_2012,ben-tal_robust_2009}} \end{spacing} }  }
      child [grow=170, level distance=125] 
        {node [concept] (other) { \begin{spacing}{0} Other\\Approaches\\{\tiny \citep{esfahani_data-driven_2015,jiang_data-driven_2016,shapiro_distributionally_2017} }  \end{spacing} }};
  \end{scope}

  \path[mindmap,concept color=black!80,text=white]
    node[concept] (chance-constraint) {
    {\Large Chance-Constrained Optimization \par} {\tiny \citep{charnes_cost_1958} }
    }
    [clockwise from=0]
    child[concept color=red1, grow=45, level distance=130]
      { node[concept] {\begin{spacing}{0} {\normalsize RO-based Methods} {\tiny \citep{ben-tal_robust_2009,bertsimas_theory_2011}} \end{spacing} } 
      [clockwise from=0]
      child[concept color=red2, grow=15]
        { node[concept] (convex-approx) { \begin{spacing}{0} Convex\\Approximation\\{\tiny\citep{nemirovski_convex_2006}} \end{spacing} }
          child[concept color=red3, grow=-45] 
            { node[concept] (cvar-chen-2010) { \begin{spacing}{0} CVaR {\tiny\citep{chen_cvar_2010}} \end{spacing} } }
          child[concept color=red3, grow=0]
            { node[concept] { \begin{spacing}{0} {\tiny\citep{chen_goal-driven_2009} } \end{spacing} } }  
          child[concept color=red3, grow=45]
            { node[concept] { \begin{spacing}{0} {\tiny\citep{rockafellar_optimization_2000,rockafellar_conditional_2002} } \end{spacing} } }  
        }
      child[concept color=red2, grow=-45]
        { node[concept] (robust-counterpart) { \begin{spacing}{0} Robust Counterpart {\tiny \citep{ben-tal_robust_2009}} \end{spacing} } 
          child[concept color=red3, grow=-15]
            { node[concept] { \begin{spacing}{0} Data-driven RO\\{\tiny\citep{bertsimas_data-driven_2018}} \end{spacing} } }
          child[concept color=red3, grow=-60]
            { node[concept] { \begin{spacing}{0} {\tiny\citep{bertsimas_price_2004} } \end{spacing} } }
          child[concept color=red3, grow=-150]
            { node[concept] { \begin{spacing}{0} {\tiny\citep{ben-tal_robust_1999} } \end{spacing} } }  
        }   
    } 
    child[concept color=blue1, grow=-45, level distance=135] {
      node[concept] { \begin{spacing}{0} {\normalsize Scenario Approach} {\tiny \citep{calafiore_uncertain_2005,campi_introduction_2018} } \end{spacing} }
      [clockwise from=0]
      child[concept color=blue2, grow=10] { node[concept] {
      \begin{spacing}{0}
      Prior\\
      {\tiny \citep{campi_exact_2008} \\\citep{calafiore_random_2010} }   
      \end{spacing}
      }
      child[concept color=blue3, grow=130]{ node[concept] (margellos-road-2014) { {\tiny\citep{margellos_road_2014}} } }
      child[concept color=blue3, grow=40] { node[concept] { {\tiny\citep{zhang_sample_2015} } } }
      child[concept color=blue3, grow=-5] { node[concept] { {\tiny\citep{schildbach_randomized_2013} } } }
      child[concept color=blue3, grow=-50] { node[concept] { {\tiny\citep{ramponi_consistency_2018} } } }
      }
      child[concept color=blue2, grow=-45] { node[concept] {
      \begin{spacing}{0}
        Posterior {\tiny\citep{campi_wait-and-judge_2016}}
      \end{spacing}
      }
        child[concept color=blue3, grow=-30] { node[concept] { {\tiny\citep{campi_general_2018} } } }
      }
      child[concept color=blue2, grow=-120] { node[concept] (sampling-and-discarding) {
      \begin{spacing}{0}
        Sampling \& Discarding {\tiny \citep{campi_sampling-and-discarding_2011} } 
      \end{spacing}
      } }
    }
    child[concept color=green1, grow=-135, level distance=135]
      {node[concept] (saa) {\begin{spacing}{0} Sample\\Average\\Approximation\\{\tiny \citep{sen_relaxations_1992,ahmed_solving_2008} }\end{spacing}}
      [clockwise from=-60]
      child[concept color=green2, grow=-60]
        { node[concept] (general-theory) { \begin{spacing}{0} General Theory \\{\tiny\citep{luedtke_sample_2008,pagnoncelli_sample_2009}} \end{spacing} } }
      child[concept color=green2, grow=-135]
        { node[concept] (sth) { \begin{spacing}{0} Computation \end{spacing} } 
          child[concept color=green3, grow=160]
            { node[concept] { \begin{spacing}{0} {\tiny \cite{tanner_iis_2010,kucukyavuz_mixing_2012} } \end{spacing} }}
          child[concept color=green3, grow=-150]
            { node[concept] { \begin{spacing}{0} {\tiny \cite{zeng_chance_2017,zeng_solving_2014}} \end{spacing} }}    
        }
      child[concept color=green2, grow=170]
        { node[concept] { \begin{spacing}{0} Special Cases \end{spacing} } 
          child[concept color=green3, grow=70]
            { node[concept] (ahmed-relaxations-2018) { \begin{spacing}{0} {\tiny \citep{ahmed_relaxations_2018} } \end{spacing} }}
          child[concept color=green3, grow=150]
            { node[concept] (luedtke-integer-2010) { \begin{spacing}{0} {\tiny \citep{luedtke_integer_2010}} \end{spacing} }}  
        }
    }
    child[concept color=brown1, grow=135, level distance=130] {
      node[concept] {\begin{spacing}{0} {\normalsize Basic Properties} {\tiny \citep{shapiro_lectures_2009} } \end{spacing} }
      [clockwise from=0]
      child[concept color=brown2, grow=-135] { node[concept]
      { \begin{spacing}{0} Special Cases \end{spacing} } 
        child[concept color=brown3, grow=0] { node[concept]
          { \begin{spacing}{0} Gaussian {\tiny \citep{kataoka_stochastic_1963}} \end{spacing} } }
        child[concept color=brown3, grow=180] { node[concept] (log-concave)
          { \begin{spacing}{0} Log-concave {\tiny \citep{prekopa_stochastic_1995}} \end{spacing} } }
        child[concept color=brown3, grow=-45] { node[concept] (finite-distribution)
          { \begin{spacing}{0} Finite Distribution {\tiny \citep{prekopa_stochastic_1995}} \end{spacing} } } 
        child[concept color=brown3, grow=-140] { node[concept] (separable-cc)
          { \begin{spacing}{0} Separable Chance Constraint {\tiny \citep{prekopa_stochastic_1995}} \end{spacing} } }       
      }
      child[concept color=brown2, grow=180]
        { node[concept] (convexity) {\begin{spacing}{0} Convexity\end{spacing} }
          child[concept color=brown3, grow=120]
            { node[concept] (calafiore-distributionally-2006) {\begin{spacing}{0} {\tiny \citep{calafiore_distributionally_2006} } \end{spacing} }}
          child[concept color=brown3, grow=180]
            { node[concept] {\begin{spacing}{0} {\tiny \citep{lagoa_convexity_1999} } \end{spacing} }}              
        }
      };
  \begin{pgfonlayer}{background}
    \draw (general-theory) to[circle connection bar switch color=from (green2) to (blue2)] (sampling-and-discarding);
    \draw (convexity) to[circle connection bar switch color=from (brown2) to (brown3)] (log-concave);
    \draw (separable-cc) to[circle connection bar switch color=from (brown3) to (green3)] (luedtke-integer-2010);
    \draw (separable-cc) to[circle connection bar switch color=from (brown3) to (green3)] (ahmed-relaxations-2018);
    \draw (finite-distribution) to[circle connection bar switch color=from (brown3) to (green3)] (ahmed-relaxations-2018);
    \draw (margellos-road-2014) to[circle connection bar switch color=from (blue3) to (red2)] (robust-counterpart);
    \draw (safe-approx) to[circle connection bar switch color=from (lightgray) to (red2)] (convex-approx);
    \draw (ambiguous-cc) to[circle connection bar switch color=from (lightgray) to (black!80)] (chance-constraint); 
    \draw (cvar-chen-2010) to[circle connection bar switch color=from (red3) to (red2)] (robust-counterpart);
  \end{pgfonlayer}
\end{tikzpicture}



\end{adjustbox}
\caption{A Schematic Overview of Existing Methods and Algorithms to Solve Chance-constrained Optimization Problems}
\label{fig:overview_cco}
\end{figure*}

A schematic overview of solutions to (CCO) and their relationships are presented in Figure \ref{fig:overview_cco}. Akin methods are plotted in similar colors, and links among two circles indicate the connection of the two methods.
The tree-like structure of Figure \ref{fig:overview_cco} illustrates the hierarchical relationship of the reviewed methods.  Key references of each method are also provided. The root node of Figure \ref{fig:overview_cco} is the ``ambiguous chance constraint'' or distributionally robust optimization (DRO), which is the parent node of ``chance-constrained optimization''. This indicates that DRO contains CCO as a special case. Similarly, for example, node ``scenario approach'' has three child nodes ``prior'', ``posterior'' and ``sampling and discarding'', this indicates the scenario approach has three major variations.

As shown in Figure \ref{fig:overview_cco}, CCO is a special case of ambiguous chance constraints where the set of distributions $\mathcal{P}$ is a singleton (Section \ref{sub:ambiguous_chance_constraints}). Therefore methods to solve ambiguous chance constraints can be applied on chance constraints as well. The methods and algorithms to solve CCO are the main focus of this paper, we will briefly mention the connection if some methods are related with ambiguous chance constraints.

Figure \ref{fig:overview_cco} also outlines the first half of this paper, which dedicates to a review and tutorial on chance-constrained optimization. We summarize key results on the basic properties (Section \ref{sec:properties_chance_constrained_optimization}), three main approaches to solving chance-constrained optimization problems, scenario approach (Section \ref{sec:scenario_approach}), sample average approximation (Section \ref{sec:sample_average_approximation}) and robust optimization (RO) based methods (Section \ref{sec:robust_optimization}).


\section{Scenario Approach} 
\label{sec:scenario_approach}
\subsection{Introduction to Scenario Approach} 
\label{sub:introduction_to_scenario_approach}
Scenario approach utilizes a dataset with $N$ scenarios $\{\xi^i\}_{i=1}^N$ to approximate the chance-constrained program (\ref{form:cco_joint}) and obtains the following \emph{scenario problem} $\text{(SP)}_N$:
\begin{subequations}
\label{form:scenario_problem}
\begin{align}
  \text{(SP)$_N$:}~\min_{x \in \mathcal{X}}~&  c^\intercal x  \\
  \text{s.t.}~&  f(x, \xi^1) \le 0, \cdots, f(x, \xi^N) \le 0  
\end{align}
\end{subequations}
$\text{SP}_N$ seeks the optimal solution $x_N^\ast$ which is feasible for all $N$ scenarios. The scenario approach is a very simple yet powerful method. The most attractive feature of the scenario approach is its generality. It requires nothing except the convexity of constraints $f(x,\xi)$ and $\mathcal{X}$. It is purely data-driven and makes no assumption on the underlying distribution.
\begin{rem}
\label{rem:SP_N_is_a_random_program}
$\text{SP}_N$ is a random program. Both its optimal objective value $o^\ast_N$ and optimal solution $x_N^\ast$ depend on the random samples $\{\xi^i\}_{i=1}^{N}$, therefore they are random variables. In consequence, $\mathbb{V}(x_N^\ast)$ is also a random variable. Let $\mathcal{N} := \{1,2,\cdots,N\}$ denote the index set of scenarios. The optimal objective value of $\text{SP}_N$ is denoted by $o^\ast(\mathcal{N})$ to emphasize its dependence on the random samples.
\end{rem}
Theoretical results of the scenario approach are built upon the following assumption in addition to Assumptions \ref{assu:distribution}, \ref{assu:convexity} and \ref{assu:iid_samples}.
\begin{assumption}[Feasibility and Uniqueness \citep{campi_exact_2008}]
\label{asmp:feasibility_uniqueness}
Every scenario problem $\text{(SP)}_N$ is feasible, and its feasibility region has a non-empty interior. Moreover, the optimal solution $x_N^\ast$ of $\text{(SP)}_N$ exists and is unique.
\end{assumption}
If there exist multiple optimal solutions, the tie-break rules in \citep{calafiore_uncertain_2005} can be applied to obtain a unique solution.

\begin{rem}[Sample Complexity $N$]
\label{rem:sample_complexity_question}
We first provide some intuition on the scenario approach. When solving $\text{(SP)}_N$ with a very large number of scenarios, the solution $x_N^\ast$ will be robust to almost every realization of $\xi$, thus the violation probability goes to zero. Although $x_N^\ast$ is a feasible solution to (CCO) as $N\rightarrow +\infty$, it is overly conservative because $\mathbb{V}(x^\ast) \approx 0 \ll \epsilon$. On the other hand, using too few scenarios for $\text{SP}_N$ might result in infeasible solutions $x_N^\ast$ to (CCO). Notice that $N$ is the only tuning parameter in the scenario approach, the most important question in the scenario approach theory is: \emph{what is the right sample complexity $N$?} Namely, what is the smallest $N$ such that $\mathbb{V}(x_N^\ast) \le \epsilon$ (with high probability)? Rigorous answers to the sample complexity question are built upon the structural properties of $\text{SP}_N$.
\end{rem}

\subsection{Structural Properties of the Scenario Problem} 
\label{sub:structural_properties_of_the_scenario_problem}
Among $N$ scenarios in the dataset $\{\xi^i\}_{i=1}^{N}$, there are some important scenarios having direct impacts on the optimal solution $x_N^\ast$.
\begin{defn}[Support Scenario \citep{calafiore_uncertain_2005}]
Scenario $\xi^i$ is a \emph{support scenario} for $\text{(SP)}_N$ if its removal changes the solution of $\text{(SP)}_N$. The set of \emph{support scenarios} of $(\text{SP}_N)$ is denoted by $\mathcal{S}$.
\end{defn}
\begin{thm}[\citep{calafiore_uncertain_2005,calafiore_random_2010}]
\label{thm:max_num_support_scenario}
Under Assumption \ref{assu:convexity}, the number of support scenarios in $\text{SP}_N$ is at most $n$, i.e. $|\mathcal{S}| \le n$.
\end{thm}
Theorem \ref{thm:max_num_support_scenario} is built upon Helly's theorem and Radon's theorem \citep{rockafellar_convex_2015} in convex analysis. For non-convex problems, the number of support scenarios could be greater than the number of decision variables $n$. An example for non-convex problems is provided in \citep{campi_general_2018}.
\begin{defn}[Fully-supported Problem \citep{campi_exact_2008}]
A scenario problem $\text{SP}_N$ with $N \ge n$ is \emph{fully-supported} if the number of support scenarios is exactly $n$. Scenario problems with $|\mathcal{S}| < n$ are referred as \emph{non-fully-supported} problems.
\end{defn}

\begin{defn}[Non-degenerate Problem \citep{campi_exact_2008,calafiore_random_2010}]
Problem $\text{SP}_N$ is said to be \emph{non-degenerate}, if $o^\ast(\mathcal{N}) = o^\ast( \mathcal{S} )$. In other words, $\text{SP}_N$ is \emph{non-degenerate} if the solution of $\text{(SP)}_N$ with all scenarios in place coincides with the solution to the program with only the support scenarios are kept.
\end{defn}

\subsection{A-priori Feasibility Guarantees} 
\label{sub:a_priori_scenario approach}
Obtaining a-priori feasibility guarantees on the solution $x_N^\ast$ to $\text{SP}_N$ typically involves the following three steps:
\begin{enumerate}
	\item Exploring the problem structure of $\text{SP}_N$ and obtain an upper bound $\overline{h}$ on the number of support scenarios;
	\item Choosing a good sample complexity $N(\epsilon, \beta,\overline{h})$ using Corollary \ref{cor:prior_sample_complexity_fully_supported}, Theorem \ref{thm:prior_guarantee_helly_dim} or Remark \ref{rem:sample_complexity_revisited};
	\item Solving the scenario problem $\text{SP}_N$ and obtain $x_N^\ast$ and $o_N^\ast$.
\end{enumerate}
\begin{thm}[\citep{campi_exact_2008}]
\label{thm:exact_feasibility_scenario_approach}
Under Assumption \ref{assu:distribution}, \ref{assu:convexity} and \ref{assu:iid_samples}, for a non-degenerate problem $\text{SP}_N$, it holds that
\begin{equation}
\label{eqn:exact_distribution_violation_prob}
	\mathbb{P}^N \Big( \mathbb{V}(x_N^\ast)  > \epsilon \Big) \le \sum_{i=1}^{n-1} \binom{N}{i} \epsilon^i (1- \epsilon )^{N-i}.
\end{equation}
The probability $\mathbb{P}^N$ is taken with respect to $N$ random samples $\{\xi^i\}_{i=1}^N$, and the inequality is tight for fully-supported problems.
\end{thm}
As mentioned in Remark \ref{rem:SP_N_is_a_random_program}, $\mathbb{V}(x_N^\ast)$ is a random variable, its randomness comes from drawing scenarios $\{\xi^i\}_{i=1}^N$. For fully-supported problems, Theorem \ref{thm:exact_feasibility_scenario_approach} shows the \emph{exact} probability distribution of the violation probability $\mathbb{V}(x_N^\ast)$, i.e.
\begin{equation}
	\mathbb{P}^N \Big( \mathbb{V}(x_N^\ast)  > \epsilon \Big) = \sum_{i=1}^{n-1} \binom{N}{i} \epsilon^i (1- \epsilon )^{N-i},
\end{equation}
the tail of a binomial distribution. We could use Theorem \ref{thm:exact_feasibility_scenario_approach} to answer the sample complexity question in Remark \ref{rem:sample_complexity_question}.
\begin{cor}[\citep{campi_exact_2008}]
\label{cor:prior_sample_complexity_fully_supported}
Given a violation probability $\epsilon \in (0,1)$ and a confidence parameter $\beta \in (0,1)$, if we choose the number of scenarios $N$ (the smallest such $N$ is denoted by $N_{2008}$) such that 
\begin{equation}
  \sum_{i=0}^{n-1} \binom{N}{i} \epsilon^i (1- \epsilon)^{N-i} \le \beta
\end{equation}
Let $x_N^\ast$ denote the optimal solution to $\text{SP}_N$, it holds that
\begin{equation}
  \mathbb{P}^N\Big ( \mathbb{V}(x_N^\ast) \le \epsilon \Big) \ge 1 - \beta
\end{equation}	
In other words, the optimal solution $x_N^\ast$ is a feasible solution to (CCO) with probability at least $1- \beta$.
\end{cor}
For fully-supported problems, $N_{2008}$ is the tightest upper bound on sample complexity, which cannot be improved. For non-fully supported problems, it turns out $N_{2008}$ can be further tightened. An improved sample complexity bound is provided in Theorem \ref{thm:prior_guarantee_helly_dim} based on the definition of Helly's dimension.
\begin{defn}[Helly's Dimension \citep{calafiore_random_2010}]
\emph{Helly's dimension} of $\text{SP}_N$ is the smallest integer $h$ such that 
\begin{equation*}
  \esssup_{\xi \in \Xi^N} | \mathcal{S}(\xi) | \le h
\end{equation*}
holds for any finite $N \ge 1$. The essential supremum is denoted by $\esssup$. We emphasize the dependence of support scenarios $\mathcal{S}$ on $\xi$ by $\mathcal{S}(\xi)$.
\end{defn}
\begin{thm}[\citep{calafiore_random_2010}]
\label{thm:prior_guarantee_helly_dim}
Let $h$ denote the Helly's dimension for $\text{SP}_N$, under Assumption \ref{assu:distribution}, \ref{assu:convexity} and \ref{assu:iid_samples}, for a non-degenerate problem $\text{SP}_N$, it holds that
\begin{eqnarray}
\label{eqn:prior_guarantee_helly_dim}
  \mathbb{P}^N \big( \mathbb{V}(x_N^\ast) > \epsilon \big) \le \sum_{i=0}^{h-1} \binom{N}{i} \epsilon^i (1- \epsilon)^{N-i}
\end{eqnarray}
Equivalently, for a fixed confidence parameter $\beta \in (0,1)$, if the sample complexity $N$ satisfies
\begin{eqnarray}
\label{eqn:prior_guarantee_binomial_tail}
  \sum_{i=0}^{h-1} \binom{N}{i} \epsilon^i (1- \epsilon)^{N-i} \le \beta
\end{eqnarray}
then the following probabilistic guarantee holds
\begin{eqnarray}
  \mathbb{P}^N \big( \mathbb{V}(x_N^\ast) > \epsilon \big) \le \beta
\end{eqnarray}
\end{thm}
The only difference between Theorem \ref{thm:prior_guarantee_helly_dim} and Theorem \ref{thm:exact_feasibility_scenario_approach} (and Corollary \ref{cor:prior_sample_complexity_fully_supported}) is replacing $n$ with Helly's dimension $h$ in (\ref{eqn:prior_guarantee_helly_dim}) and (\ref{eqn:prior_guarantee_binomial_tail}). Unfortunately, Helly's dimension is often difficult to calculate, while finding upper bounds $\overline{h}$ on Helly's dimension is usually a much easier task. Similarly we can replace $h$ by $\overline{h}$ in (\ref{eqn:prior_guarantee_helly_dim}) and (\ref{eqn:prior_guarantee_binomial_tail}), the same theoretical guarantees still hold because of the monotonicity of (\ref{eqn:prior_guarantee_helly_dim}) and (\ref{eqn:prior_guarantee_binomial_tail}) in $N$ and $h$. The support-rank defined in \cite{schildbach_randomized_2013} is an upper bound on Helly's dimension, some other upper bounds can be obtained by exploiting the structural properties of the problem, e.g. \cite{zhang_sample_2015}.
\begin{rem}[Sample Complexity Revisited]
\label{rem:sample_complexity_revisited}
A binary search type algorithm could be used to find $N_{2008}$. And a looser but handy upper bound is provided in \citep{campi_scenario_2009}:
\begin{equation}
\label{eqn:N_2009}
	N_{2009} := \frac{2}{\epsilon}\Big(\ln(\frac{1}{\beta}) + n \Big)
\end{equation}
Notice $n$ in (\ref{eqn:N_2009}) can be replaced by $h$ or $\overline{h}$.
\end{rem}

\subsection{A-posteriori Feasibility Guarantees} 
\label{sub:a_posterior_Guarantees}
When the desired violation probability $\epsilon$ is very small, the sample complexity of the a-priori guarantees grows with $1/\epsilon$ (Remark \ref{rem:sample_complexity_revisited}) and could be prohibitive. In other words, the a-priori approach is only suitable for the case where a \emph{sufficient} amount of scenarios is always available. In many real-world applications (e.g. medical experiments, tests conducted by NASA), however, the amount of data is quite limited, and it could take months or cost a fortune to obtain a data point (experiment).
Because of the limitation on the data availability, one of the most fundamental problem in data-driven decision making (e.g. system identification, quantitative finance) is to come up with good decisions or estimates with a moderate or even small amount of data. To overcome this, the scenario approach is extended towards a-posteriori feasibility guarantees. 

Similar with the \emph{a-priori} guarantees, obtaining a-posteriori guarantees typically requires taking the following three steps:
\begin{enumerate}
	\item given dataset $\{\xi^i\}_{i=1}^N$, solve the corresponding scenario problem $\text{SP}_N$ and obtain $x_N^*$;
	\item find support scenarios in $\{\xi^i\}_{i=1}^N$, whose number is denoted as $s_N^*$;
	\item calculate the posterior violation probability $\epsilon(\beta, s_N^*,N)$ using Theorem \ref{thm:posterior_guarantee}.
\end{enumerate} 
If the resulting violation probability $\epsilon(\beta, s_N^*,N)$ is greater than the acceptable level $\epsilon$, we could repeat this process with more scenarios until reaching $\epsilon(\beta, s_N^*,N) \le \epsilon$. If the number of available scenarios is limited, then it might be impossible to obtain a solution $x_N^*$ such that $\mathbb{V}(x_N^*) \le \epsilon$.
\begin{thm}[Wait-and-Judge \citep{campi_wait-and-judge_2016}]
\label{thm:posterior_guarantee}
Given $\beta \in (0,1)$, for any $k=0,1,\cdots,n$, the polynomial equation in variable $t$
\begin{equation}
	\frac{\beta}{N+1} \sum_{i=k}^N \binom{i}{k} t^{i-k} - \binom{N}{k} t^{N-k} = 0
\end{equation}
has exactly one solution $\epsilon(k)$ in the interval $(0,1)$. Under Assumption \ref{assu:distribution}, \ref{assu:convexity} and \ref{assu:iid_samples}, for a non-degenerate problem, it holds that
\begin{equation}
	\mathbb{P}^N(\mathbb{V}(x_N^\ast) \ge \epsilon(s_N^\ast)) \le \beta
\end{equation}
\end{thm}
Theorem \ref{thm:posterior_guarantee} is particularly useful in the following cases: (i) the problem is not fully-support thus difficult to calculate \emph{a-priori} bounds on number of support scenarios; or (ii) only a moderate or small amount of data points is available, it is difficult to meet the sample complexity from the a-priori guarantees. 




Given a candidate solution $x^\diamond$, the most straightforward method is to approximate $\mathbb{V}(x^\diamond)$ is by the empirical estimation $\hat{\epsilon}$ through Monte-Carlo simulation with $\hat{N}$ samples, i.e.
\begin{equation}
\label{eqn:empirical_estimate_violation_prob}
	\hat{\epsilon} = \frac{1}{\hat{N}} \sum_{i=1}^{\hat{N}} \mathbbm{1}_{f(x^\diamond, \xi^i) > 0} = \frac{\hat{V}}{\hat{N}}
\end{equation}
where $\hat{V} := \sum_{i=1}^{\hat{N}} \mathbbm{1}_{f(x^\diamond, \xi^i) > 0}$ is the total number of scenarios in which $x^\ast_N$ is infeasible. Although (\ref{eqn:empirical_estimate_violation_prob}) only involves $f(x^\diamond, \xi^i) > 0$ which is easy to calculate, it might require an astronomical number $\hat{N}$ to have accurate estimation $\hat{\epsilon}$ because of (D1). \citep{nemirovski_convex_2006} shows a method to bound $\mathbb{V}(x^\diamond)$ from above using a dataset of a moderate size $\hat{N}$.
\begin{prop}[\citep{nemirovski_convex_2006}]
\label{prop:check_feasibility_scenario_approach}
Given a candidate solution $x^\diamond$ and $\hat{N}$ samples, let $\hat{V} := \sum_{i=1}^{\hat{N}} \mathbbm{1}_{f(x^\diamond, \xi^i) > 0}$ and $1-\rho$ be the confidence parameter.
\begin{equation}
\label{eqn:upper_bound_violation_prob}
	\overline{\epsilon} := \max_{\gamma \in [0,1]} \{\gamma: \sum_{i=0}^{\hat{V}} \binom{\hat{N}}{i} \gamma^i (1- \gamma)^{\hat{N}-i} \ge \rho  \}
\end{equation}
After finding an upper bound $\overline{\epsilon}$, so that if $\overline{\epsilon} \le \epsilon$, we may be sure that $\mathbb{P}(\mathbb{V}(x^\diamond) \le \epsilon) \ge 1 - \rho$. 
\end{prop}
\begin{rem}
Proposition \ref{prop:check_feasibility_scenario_approach} is closely related with scenario approach but with one fundamental difference. Theorem \ref{thm:posterior_guarantee} holds only for solution from scenario approach, while Proposition \ref{prop:check_feasibility_scenario_approach} can evaluate solutions from other methods.
\end{rem}

\subsection{Optimality Guarantees of Scenario Approach} 
\label{sub:constructing_lower_bounds_via_scenario_approach}
Scenario approach together with order statistics can be used to construct lower bounds $\underline{o}$ on $o^\star$ of (CCO).
\begin{prop}[\citep{nemirovski_convex_2006}]
\label{prop:lower_bound_cco_framework}
Let $\{\xi^{i,j}\}_{i=1}^N$ ($j=1,2,\cdots,M$) be $M$ \emph{independent} datasets of size $N$. For the $j$ th dataset, we solve the associated scenario problem $\text{SP}_N$ and calculate the optimal value $o_j^\ast$ ($j=1,2,\cdots,M$). Without loss of generality, we assume that $o_1^\ast \le o_2^\ast \le \cdots \le o_M^\ast$.

Given $\delta \in (0,1)$, let us choose positive integers $M$,$N$,$L$ in such a way that 
\begin{equation}
\label{eqn:order_statistics_lower_bound}
\sum_{i=0}^{L-1} \binom{M}{i} (1- \epsilon)^{Ni}[1-(1- \epsilon)^N]^{M-i} \le \delta
\end{equation}
then with probability of at least $1- \delta$, the random quantity $o_L^*$ gives a lower bound for the true optimal value $x^\star$.
\end{prop}
\citep{pagnoncelli_sample_2009} shows that appropriate $N$ should be the order of $O(1/\epsilon)$ as $[1-(1- \epsilon)^N]^M \approx (1- \exp(-\epsilon N))^M$. Typically we choose proper values for $N$ and $M$ first, then find out the largest positive integer $L$ that (\ref{eqn:order_statistics_lower_bound}) holds true.

Proposition \ref{prop:lower_bound_cco_framework} turns out to be a general framework to construct lower bounds on (CCO). \citep{pagnoncelli_sample_2009} extends the framework towards generating bounds using sample average approximation, which is introduced in Section \ref{sub:optimality_guarantees_of_sample_average_approximation}.




\section{Sample Average Approximation} 
\label{sec:sample_average_approximation}

\subsection{Introduction to Sample Average Approximation} 
\label{sub:introduction_to_sample_average_approximation}
The idea of using sample average approximation to handle chance constraints first appeared in \citep{sen_relaxations_1992} and was subsequently improved with rigorous theoretical results in \citep{luedtke_sample_2008}.

Let $\overline{f}(x,\xi):=\max\big\{f_1(x,\xi),\cdots,f_m(x,\xi) \big\}$, then (CCO) is equivalent to $\min_{x \in \mathcal{X}} c^\intercal x,~\text{s.t.}~\mathbb{P}(\overline{f}(x,\xi) \le 0) \ge 1 - \epsilon$. Sample Average Approximation (SAA) approximates the true distribution of the random variable $\overline{f}(x,\xi)$ using the empirical distribution from $N$ samples $\{\xi^i\}_{i=1}^N$, i.e. $\mathbb{P}(\overline{f}(x,\xi) \le 0)$ is approximated by $\frac{1}{N} \sum_{i=1}^N \mathbbm{1}_{\overline{f}(x, \xi^i) \le 0}$.
\begin{subequations}
\begin{align}
\text{(SAA):} \min_{x\in \mathcal{X}}~& c^\intercal x \\
\text{s.t.}~& \frac{1}{N} \sum_{i=1}^N \mathbbm{1}_{\overline{f}(x, \xi^i) > 0}  \le \varepsilon
\end{align}
\end{subequations}
(SAA) is also a chance constrained optimization problem, but with two major differences from (CCO): (i) (SAA) is based on the empirical (discrete) distribution from the true distribution of $\xi$ as in (CCO); (ii) (SAA) has the violation probability $\varepsilon$ instead of $\epsilon$ in (CCO).

There are two critical questions to be addressed about (SAA). What is the connection of solutions of (SSA) with that of (CCO)? How to solve (SAA)? We first answer the second question in Section \ref{sub:solving_sample_average_approximation}, then present the theoretical results of connecting (SAA) with (CCO).

\subsection{Solving Sample Average Approximation} 
\label{sub:solving_sample_average_approximation}
(SAA) can be reformulated as a mixed integer program (MIP) by introducing variables $z \in \{0,1\}^N$ \citep{ruszczynski_probabilistic_2002,luedtke_sample_2008}. Binary variable $z_i$ is an indicator if $\overline{f}(x,\xi) \le 0$ is being violated in sample $i$, i.e.
\begin{equation}
\label{eqn:indicator_z}
z_i = \mathbbm{1}_{\overline{f}(x,\xi^i) > 0}
\end{equation}
(\ref{eqn:indicator_z}) can be equivalently written as $\overline{f}(x,\xi^i) \le M z_i$ with a sufficiently large coefficient $M \in \mathbf{R}_+$. Since $\overline{f}(x,\xi^i)$ is the maximum over $m$ functions $\{f_j(x,\xi^i)\}_{j=1}^m$, $\overline{f}(x,\xi^i) \le M z_i$ implies $f_j(x,\xi^i) \le M z_i, j=1,2,\cdots,m$. Then (SAA) is equivalent to (\ref{form:saa_mip}), in which $\mathbf{1}_m$ is an all one vector with size $m$.
\begin{subequations}
\label{form:saa_mip}
\begin{align}
	\min_{x, z}~&  c^\intercal x \\
  \text{s.t.}~& f(x, \xi^1) - M z_1 \mathbf{1}_m \le 0\\
  & \vdots \nonumber \\
  & f(x, \xi^N) - M z_N  \mathbf{1}_m \le 0 \\
  & \frac{1}{N} \sum_{i=1}^N z_i \le \varepsilon  \label{form:saa_mip_cc} \\
  & x \in \mathcal{X}, z_i \in \{0,1\}, i = 1,2,\cdots,N
\end{align}
\end{subequations}
(\ref{form:saa_mip}) is equivalent to (SAA) for general function $f(x,\xi)$, but formulations with big-M are typically weak formulations. Introducing big coefficients $M$ might cause numerical issues as well. Stronger formulations of (SAA) are possible by exploiting the structural features of $f(x,\xi)$. A good example is chance-constrained linear program with separable probabilistic constraints: $\min_{x \in \mathcal{X}} c^\intercal x~\text{s.t.} \mathbb{P}(T x \ge \xi ) \ge 1- \epsilon$, with a constant matrix $T \in \mathbf{R}^{d\times n}$. By introducing auxiliary variables $v$, an equivalent but stronger formulation without big M is (\ref{form:cclp_saa_mip_noM}) \citep{luedtke_integer_2010}.
\begin{subequations}
\label{form:cclp_saa_mip_noM}
\begin{align}
\min_{x \in \mathcal{X}}~& c^\intercal x \\
\text{s.t.}~& Tx = v \\
& v + \xi_i z_i \ge \xi_i, i = 1,2,\cdots,N\\
& \frac{1}{N} \sum_{i=1}^N z_i \le \varepsilon \\
& z_i \in \{0,1\}, i= 1,2,\cdots,N 
\end{align}
\end{subequations}
Various strong formulations for (SAA) can be found in \citep{luedtke_integer_2010} and references therein. (\ref{form:saa_mip}) and (\ref{form:cclp_saa_mip_noM}) are mixed integer programs, some well-known techniques from integer programming theory can speed up the process of solving (SAA), e.g. adding cuts \citep{tanner_iis_2010,luedtke_integer_2010,kucukyavuz_mixing_2012} and decompositions \citep{zeng_solving_2014,zeng_chance_2017}.

\subsection{Feasibility Guarantees of SAA} 
Various feasibility guarantees of (SAA) are proved in \citep{luedtke_sample_2008,pagnoncelli_sample_2009}, e.g. the asymptotic behavior of (SAA), when $\mathcal{X}$ is finite, the case of separable chance constraints (\ref{form:cco_separable}), and when $f(x,\xi)$ is Lipschitz continuous. In this section, we only present the Lipschitz case, which could be useful for many engineering applications.
\begin{assumption}
\label{assu:lipschitz}
There exists $L >0$ such that 
\begin{equation}
	|\overline{f}(x,\xi) - \overline{f}(x^\prime,\xi)| \le L \|x-x^\prime\|_\infty, ~\forall x,x^\prime \in \mathcal{X}~\text{and}~\forall \xi \in \Xi.
\end{equation}
\end{assumption}
\begin{thm}[\citep{luedtke_sample_2008}]
\label{thm:lipschitz-saa}
Suppose $\mathcal{X}$ is bounded with diameter $D$ and $\overline{f}(x,\xi)$ is $L$-Lipschitz for any $\xi \in \Xi$ (Assumption \ref{assu:lipschitz}). Let $\varepsilon \in [0, \epsilon), \theta \in (0,\epsilon - \varepsilon)$ and $\gamma > 0$. Then 
\begin{equation}
	\mathbb{P}(\mathcal{F}_{\varepsilon, \gamma}^N \subseteq \mathcal{F}_\epsilon ) \ge 
	1 - \ceil{\frac{1}{\theta}} \ceil{\frac{2LD}{\gamma}}^n \exp(-2N(\epsilon - \varepsilon - \theta)^2)
\end{equation}
where the feasible region of (SAA) is defined as
\begin{equation}
	\mathcal{F}_{\varepsilon, \gamma}^N := \{ x\in \mathcal{X}: \frac{1}{N} \sum_{i=1}^N \mathbbm{1}_{\overline{f}(x,\xi) + \gamma \le 0} \ge  1 - \varepsilon \}.
\end{equation}
\end{thm}
For fixed $\epsilon$ and $\varepsilon$, if we choose $\theta = (\epsilon - \varepsilon)/2$ and a small number $\gamma > 0$, then Theorem \ref{thm:lipschitz-saa} suggests that using 
\begin{equation}
	N \ge \frac{2}{(\epsilon - \varepsilon)^2} \big[ \ln(\frac{1}{\beta}) + n \ln(\ceil{\frac{2LD}{\gamma}}) + \ln(\ceil{\frac{2}{\epsilon - \varepsilon}}) \big]
\end{equation}
number of samples, solutions of (SAA) is feasible to (CCO) with high probability $1 - \beta$, i.e. 
 $\mathbb{P}(\mathcal{F}_{\varepsilon, \gamma}^N \subseteq \mathcal{F}_\epsilon ) \ge 1 - \beta$. 

The results in Theorem \ref{thm:lipschitz-saa} look quite similar to those of scenario approach (e.g. Remark \ref{rem:sample_complexity_revisited}). Indeed, (SAA) with $\varepsilon=0$ is exactly the same as the scenario problme $\text{SP}_N$. However, one major difference of Theorem \ref{thm:lipschitz-saa} from the scenario approach theory is that: Theorem \ref{thm:lipschitz-saa} holds for every feasible point of (SAA), i.e. $\mathcal{F}_{\varepsilon, \gamma}^N \subseteq \mathcal{F}_\epsilon$ with high probability. While the theory of the scenario approach only proves the property of the optimal solution $x_N^*$, i.e. $x_N^*$ is feasible with high probability. Other feasible solutions to $\text{SP}_N$ do not necessarily process the properties guaranteed by the scenario approach (e.g. Theorem \ref{thm:exact_feasibility_scenario_approach}).

Although Theorem \ref{thm:lipschitz-saa} provides explicit sample complexity bounds for (SAA) to obtain feasible solution, it requires some efforts to be applied, e.g. tuning parameters $(\varepsilon, \theta)$ and calculation of $L$ and $D$. \citep{campi_sampling-and-discarding_2011} provides a similar but more straightforward theoretical result.
\begin{thm}[Sampling \& Discarding \citep{campi_sampling-and-discarding_2011}]
\label{thm:sampling_and_discarding}
If we draw $N$ samples and discard any $k$ of them, then use the scenario approach with the remaining $N-k$ samples. If $N$ and $k$ satisfy
\begin{equation}
\label{eqn:sampling_discarding}
\binom{k+n-1}{k} \cdot \sum_{i=0}^{k+n-1} \binom{N}{i} \epsilon^i (1- \epsilon)^{N-i} \le \beta
\end{equation}
then $\mathbb{P}^N\Big ( \mathbb{P}_{\xi} (f(x_{N,k}^*,\xi) \le 0) \ge 1- \epsilon \Big) \ge 1 - \beta $.
\end{thm}
Given parameters $N$, $\epsilon$ and $\beta$, we find the largest $k$ that (\ref{eqn:sampling_discarding}) holds, then the solution to (SAA) with $\varepsilon = k/N$ is feasible to (CCO) with probability at least $1 - \beta$.

\subsection{Optimality Guarantees of Sample Average Approximation} 
\label{sub:optimality_guarantees_of_sample_average_approximation}
It is intuitive that if $\varepsilon > \epsilon$, then the objective values of SAA yield lower bounds to (CCO). Theorem \ref{thm:lower_bound_cco_saa} formalizes this intuition.
\begin{thm}[\citep{luedtke_sample_2008}]
\label{thm:lower_bound_cco_saa}
Let $\varepsilon > \epsilon$ and assume that (CCO) has an optimal solution. Then
\begin{equation}
\mathbb{P}\Big(\hat{o}_\varepsilon^N \le o_\epsilon^\star \Big) \ge 1 - \exp( -2N (\varepsilon - \epsilon)^2).
\end{equation}
\end{thm}
Theorem \ref{thm:lower_bound_cco_saa} directly suggests a method to construct lower bounds on (CCO).
\begin{prop}
If we choose $\varepsilon > \epsilon$ and $N \ge \frac{1}{2(\varepsilon - \epsilon)^2} \log(\frac{1}{\delta})$,
let $o_\varepsilon^\text{SAA}$ denote the objective value of (SAA), then $o_\varepsilon$ is a lower bound with probability at least $1- \delta$, i.e. $\mathbb{P}(o_{N,\varepsilon}^\ast \le o_\epsilon^\star) \ge 1 - \delta$.
\end{prop}
There is an alternative method using SAA to generate lower bounds of (CCO). \citep{luedtke_sample_2008} extends the framework in Proposition \ref{prop:lower_bound_cco_framework} towards SAA.
\begin{prop}[\citep{luedtke_sample_2008}]
\label{prop:lower_bound_cco_saa}
Take $K$ sets of $N$ independent samples $\{\xi^{i,j}\}_{i=1}^N$, ($j=1,2,\cdots,K$). For the $j$th dataset $\{\xi^{i,j}\}_{i=1}^N$, we solve the associated (SAA) problem and calculate the associated objective value $o_{N,\varepsilon,j}^\ast$ (for simplicity $o_j^\ast$ and $j=1,2,\cdots,K$). Without loss of generality, we assume that $o_1^\ast \le o_2^\ast \le \cdots \le o_K^\ast$.

Given $\delta \in (0,1)$, $\varepsilon \in [0,1)$, let us choose positive integers $N$,$L$, $K$ ($L \le K$) such that
\begin{equation}
	\sum_{i=0}^{L-1} \binom{K}{i} \big[ b(\varepsilon,\epsilon,N) \big]^i \big[ 1- b(\varepsilon,\epsilon,N) \big]^{K-i} \ge \delta
\end{equation}
where $b(\varepsilon,\epsilon,N) := \sum_{i=0}^{\floor{\varepsilon N}} \binom{N}{i} \epsilon^i (1- \epsilon)^{N-i}$.

Then $o_L^\ast$ serves as a lower bound to (CCO) with probability at least $1- \delta$.
\end{prop}

\section{Robust Optimization Related Methods} 
\label{sec:robust_optimization}
\subsection{Introduction to Robust Optimization} 
\label{sub:introduction_to_robust_optimization}
The last category of solutions to (CCO) is closely related with robust optimization (RO), its typical form is shown in (\ref{form:robust_counterpart}).
\begin{subequations}
\label{form:robust_counterpart}
\begin{align}
\text{(RC):}~\min_{x \in \mathcal{X}} ~ & c^\intercal x \\
\text{s.t.}~& f(x,\xi) \le 0,~\forall \xi \in \mathcal{U}_\epsilon \label{form:robust_counterpart_constraint}
\end{align}
\end{subequations}
(\ref{form:robust_counterpart}) finds the optimal solution which is feasible to all realizations of uncertainties in a predefined uncertainty set $\mathcal{U}_\epsilon$.
(\ref{form:robust_counterpart}) is called the Robust Counterpart (RC) of the original problem (CCO). By constructing an uncertainty set $\mathcal{U}_\epsilon$ with proper shape and size, solutions to (RC) could be suboptimal or approximate solutions to (CCO).

Designing uncertainty sets lies at the heart of robust optimization.
A good uncertainty set should meet the following two standards:
\begin{description}
  \item [(S1)] The resulting (RC) problem is computationally tractable.
  \item [(S2)] The optimal solution to (RC) is not too conservative or overly optimistic.
\end{description}
Unfortunately, (RC) of general robust convex problems (under Assumption \ref{assu:convexity}) is not always computationally tractable. For example, (RC) of a second order cone program (SOCP) with polyhedral uncertainty set is NP-Hard \citep{ben-tal_robust_1998,ben-tal_robust_2002-1,bertsimas_theory_2011}. Fortunately, robust linear programs are well-studied, and (RC) of linear programs is tractable for common choices of uncertainty sets.
Most tractability results of robust linear optimization are summarized in \citep{bertsimas_theory_2011}. For tractable formulations of general convex RO problems, various solutions can be found in \citep{bertsimas_tractable_2006,ben-tal_robust_2009}.

For simplicity, we present solutions to the following chance-constrained linear program (CCLP)
\footnote{A (seemingly) more general form of the linear chance constraint is $\mathbb{P}\Big( A(\xi) x  \le b(\xi)) \Big) \ge 1- \epsilon$, where $A(\xi)$ and $b(\xi)$ denote affine functions of $\xi$. This could be equivalently represented in  the form(\ref{form:cclp-cc}) by enforcing additional affine constraints \citep{chen_cvar_2010}}.
\begin{subequations}
\label{form:cclp}
\begin{align}
\min_{x \in \mathcal{X}}~& c^\intercal x \\
\text{s.t.}~ & \mathbb{P}_\xi\Big( x_0^i + \xi^\intercal x^i \le 0,~i=1,2,\cdots,m \Big) \ge 1- \epsilon \label{form:cclp-cc}
\end{align}
\end{subequations}
and its robust counterpart
\begin{subequations}
\label{form:rlp}
\begin{align}
\min_{x \in \mathcal{X}}~ & c^\intercal x  \\
\text{s.t.}~ & x_0^i + \xi^\intercal x^i \le 0,~\forall \xi \in \mathcal{U}_\epsilon,~i = 1,2,\cdots,m 
\end{align}
\end{subequations}
In (\ref{form:cclp}) and (\ref{form:rlp}), decision variables are $\{x_0^i, x^i\}_{i=1}^m$, where $x_0^i \in \mathbf{R}$ and $x^i \in \mathbf{R}^n$. Uncertainties are represented by $\xi \in \mathbf{R}^d$ \footnote{Notice $d=n$ in (\ref{form:cclp}) and (\ref{form:rlp}).} With a little abuse of notation, we use $x = [x_0^1,x^1,\cdots,x_0^m,x^m]^\intercal$ to represent all the decision variables.

Standard (S2) is directly connected with chance constraints, we present the connection between RO and CCO in Section \ref{sub:safe_approximation}-\ref{sub:safe_approximation_of_joint_chance_constraints}. 

\subsection{Safe Approximation} 
\label{sub:safe_approximation}
Almost every RO-related solution to (CCO) is based on the idea of safe approximation.
\begin{defn}[Safe Approximation]
\label{defn:safe_approx}
Let $\mathcal{F}$ and $\underline{\mathcal{F}}$ denote two sets of constraints. We say $\underline{\mathcal{F}}$ is a safe approximation (or inner approximation) of $\mathcal{F}$ if $\underline{\mathcal{F}} \subseteq \mathcal{F}$. 
\end{defn}
An optimization problem (SA) is called a \emph{safe approximation} of (CCO) if $\underline{\mathcal{F}} \subseteq \mathcal{F}_\epsilon$, where $\mathcal{F}_\epsilon$ represents the feasible region of (CCO) as in Definition \ref{defn:cco-feasible}.
\begin{subequations}
\label{form:cco_safe_approx}
\begin{align}
\text{(SA):}~\min_{x \in \mathcal{X}}~& c^\intercal x \\
\text{s.t.}~& x \in \underline{\mathcal{F}} 
\end{align}
\end{subequations}
$\underline{\mathcal{F}} \subseteq \mathcal{F}_\epsilon$ indicates that every solution to (SA) is \emph{feasible} to (CCO). Therefore every optimal solution to (SA) is \emph{suboptimal} to (CCO) and serves as an upper bound on (CCO). 

There are two major approaches to constructing safe approximations of the chance constraint $\mathbb{P}_\xi \big(f(x,\xi) \le 0 \big) \ge 1 - \epsilon$: (i) constructing a function $\pi(x) \ge \mathbb{P}_\xi \big(f(x,\xi) > 0 \big)$, then $\pi(x)\le \epsilon$ is a safe approximation of the chance constraint; (ii) constructing a proper uncertainty set $\mathcal{U}_\epsilon$ such that $\mathcal{F}_\epsilon \supseteq \mathcal{F}_{\mathcal{U}_\epsilon}:=\{x \in \mathbf{R}^n: f(x,\xi)\le 0,~\forall \xi \in \mathcal{U}_\epsilon\}$. Although these two approaches look quite different, Section \ref{ssub:cvar_based_convex_approximation_of_individual_chance_constraints} shows that they are closely related with each other.

We first review how to apply these two approaches to obtaining safe approximation of individual chance constraints in Section \ref{sub:safe_approximation}. Safe approximations of joint chance constraints (Section \ref{sub:safe_approximation_of_joint_chance_constraints}) are built upon the results of individual chance constraints.

\subsection{Safe Approximation of Individual Chance Constraints} 
\label{sub:safe_approximation_of_individual_chance_constraints}
RO has been quite successful in constructing safe approximations of individual chance constraints. A general form of individual chance-constrained programs is (\ref{form:cco-individual}). 
\begin{subequations}
\label{form:cco-individual}
\begin{align}
\min_{x \in \mathcal{X}}~& c^\intercal x \\
\text{s.t.}~ & \mathbb{P}_\xi\Big( f(x,\xi) \le 0 \Big) \ge 1- \epsilon \label{form:cco-individual-cc}
\end{align}
\end{subequations}
In the individual chance constraint (\ref{form:cco-individual-cc}), the inner function $f(x,\xi): \mathbf{R}^n \times \mathbf{R}^d \rightarrow \mathbf{R}^1$ is a \emph{scalar-valued} function. In Section \ref{sub:safe_approximation_of_individual_chance_constraints}, all $f(x,\xi)$ are scalar-valued functions if not specified.

Section \ref{sub:safe_approximation} outlines two different but related approaches to constructing safe approximations. The first approach is presented in Section \ref{ssub:convex_approximation}-\ref{ssub:cvar_based_convex_approximation_of_individual_chance_constraints}. The second approach is summarized in \ref{ssub:constructing_uncertainty_sets}.

\subsubsection{Convex Approximation} 
\label{ssub:convex_approximation}
Convex approximation is a general framework to build safe approximations of individual chance constraints. The idea of convex approximation first appeared in 
\citep{pinter_deterministic_1989}, then was completed in \citep{nemirovski_convex_2006}. The convex approximation framework is based on the concept of generating function.
\begin{defn}[Generating Function]
A function $\phi: \mathbf{R} \rightarrow \mathbf{R}$ is called a (one-dimensional) \emph{generating function} if it is nonnegative valued, nondecreasing, convex and satisfying the following property:
\begin{equation}
	\phi(z) > \phi(0) = 1, \forall z >0
\end{equation}
\end{defn}
The idea of convex approximation starts from the following lemma.
\begin{lem}
For a positive constant $t \in \mathbf{R}_+$ and a random variable $z \in \mathbf{R}$, it holds that
\begin{equation}
	\mathbb{E}[ \phi(t^{-1} z) ] \ge \mathbb{E}[ \mathbf{1}_{t^{-1}z \ge 0}] = \mathbb{P}_z (t^{-1} z \ge 0) = \mathbb{P}(z \ge 0)
\end{equation}
\end{lem}
Replace $z$ with $f(x,\xi)$, then $\mathbb{E}[ \phi(t^{-1} f(x,\xi)) ] \ge \mathbb{P}_\xi\Big( f(x,\xi) > 0 \Big) = \mathbb{P}_\xi\Big( t^{-1} f(x,\xi) > 0 \Big)$. In other words, $\mathbb{E}[ \phi(t^{-1} f(x,\xi)) ] \le \epsilon$ is a \emph{safe approximation} to $\mathbb{P}_\xi\Big( f(x,\xi) \le 0 \Big) \ge 1 - \epsilon$.
\begin{thm}[Convex Approximation \citep{nemirovski_convex_2006}]
\label{thm:convex_approx_safe_approx}
Let $\phi(\cdot)$ be a generating function, then (CA) is a safe approximation to (CCO).
\begin{subequations}
\label{form:convex_approx}
\begin{align}
  \text{(CA):}~\min_{x \in \mathcal{X}}~& c^\intercal x \\
  \text{s.t.}~& \inf_{t > 0} \big[ t \mathbb{E}_\xi[\phi(\frac{f(x,\xi)}{t})]- t \epsilon \big] \le 0  \label{form:convex_approx_constraint}
\end{align}
\end{subequations} 
Under Assumption \ref{assu:convexity}, (CA) is convex in $x$.
\end{thm}
\begin{rem}
We can get rid of the strict inequality $t > 0$ by approximating it using $t \ge \delta$, where $\delta$ is very small positive number (e.g. $\delta = 10^{-4}$).
Furthermore, we can show that (CA) is equivalent to (\ref{form:convex_approx_eqv}), which is convex in $(x,t)$.
\begin{subequations}
\label{form:convex_approx_eqv}
\begin{align}
  \min_{x \in \mathcal{X}, t \ge \delta}~& c^\intercal x \\
  \text{s.t.}~& t \mathbb{E}_\xi[\phi(\frac{f(x,\xi)}{t})]- t \epsilon \le 0
\end{align}
\end{subequations}
\end{rem}
Choosing a good generating functions plays a crucial role in the convex approximation framework. Choices of generating functions include: Markov bound $\phi(z) = [1+z]_+$, Chernoff bound $\phi(z) = \exp(z)$, Chebyshev bound $\phi(z) = [z+1]_+^2$ and Traditional Chebyshev bound $\phi(z) = (z+1)^2$.
The \emph{least} conservative generating function is the Markov bound $\phi(z) = [1+z]_+$ \citep{nemirovski_convex_2006,follmer_stochastic_2011}.
\begin{defn}[Conditional Value at Risk]
Conditional value at risk (CVaR) of a random variable $z$ at level $1- \epsilon$ is defined as
\begin{eqnarray}
	\CVaR(z;1- \epsilon) := \inf_{\gamma}(\gamma + \frac{1}{\epsilon}\mathbb{E}\big[[z-\gamma]_+\big])
\end{eqnarray}
\end{defn}
\begin{prop}[\citep{nemirovski_convex_2006,chen_cvar_2010}]
\label{prop:convex_approx_cvar}
(CA) with Markov bound $\phi(z) = [z+1]_+$ is equivalent to (\ref{form:convex_approx_cvar}).
\begin{subequations}
\label{form:convex_approx_cvar}
\begin{align}
\min_{x\in \mathcal{X}}~& c^\intercal x \\
\text{s.t.}~& \CVaR\big(f(x,\xi); 1- \epsilon\big) \le 0
\end{align}
\end{subequations}
\end{prop}

Section \ref{sec:chance_constrained_optimization} shows an individual chance constraint $\mathbb{P}\big(f(x,\xi) \le 0\big) \ge 1- \epsilon$ is equivalent to $\VaR(f(x,\xi); 1- \epsilon ) \le 0$. It is well-known that $\CVaR(z; 1- \epsilon) \ge \VaR(z; 1- \epsilon)$. Therefore, $\CVaR(f(x,\xi); 1- \epsilon) \le 0 $ implies $\VaR(f(x,\xi); 1- \epsilon) \le 0$. In other words, $\CVaR(f(x,\xi); 1- \epsilon) \le 0 $ is a safe approximation to both $\VaR(f(x,\xi); 1- \epsilon) \le 0$ and the chance constraint (\ref{form:cco-individual-cc}).

\begin{rem}[Sample Approximation of CVaR]
\label{rem:data_driven_cvar}
\citep{rockafellar_optimization_2000} utilizes a dataset $\{\xi^i\}_{i=1}^N$ to estimate CVaR.
\begin{subequations}
\label{form:convex_approx_cvar_data}
\begin{align}
  \min_{x \in X, t}~& c^\intercal x \\
  \text{s.t.}~&  \frac{1}{N} \sum_{i=1}^{N} [f(x,\xi^i)+t]_+ \le t \epsilon 
\end{align}
\end{subequations}
By introducing $N$ auxiliary variables, \citep{rockafellar_optimization_2000} shows that (\ref{form:convex_approx_cvar_data}) can be reformulated as a convex problem that is easy to solve. Detailed reformulation can be found in \citep{rockafellar_optimization_2000} and \ref{app-sub:sample_approximation_of_cvar}.
With a sufficient number of data points ($N$ is large enough), (\ref{form:convex_approx_cvar_data}) is a safe approximation to (CCO). However, it remains unknown about the exact requirement on the number of samples needed. The sample approximation of CVaR may not necessarily yield a safe approximation \citep{chen_cvar_2010}. 
\end{rem}
The generating function based framework in \citep{nemirovski_convex_2006} was further improved and completed in \citep{ben-tal_robust_2009,nemirovski_safe_2012}. But the methods proposed there are mainly analytical and aim at solving distributionally robust problems, which is beyond the scope of this paper. More details can be found in Figure \ref{fig:overview_cco} and references therein.

\subsubsection{CVaR-based Convex Approximation of Individual Chance Constraints} 
\label{ssub:cvar_based_convex_approximation_of_individual_chance_constraints}
As pointed out in \citep{nemirovski_convex_2006}, calculating CVaR is computationally intractable.
In order to obtain tractable forms of the CVaR-based convex approximation, one approach is the sample approximation in Remark \ref{rem:data_driven_cvar}. An alternative approach is to bound the CVaR function from above, e.g. finding a function $\pi(x) \ge \CVaR(f(x,\xi); 1 - \epsilon)$, then $\pi(x) \le 0$ is a safe approximation to both $\CVaR(f(x,\xi); 1- \epsilon) \le 0$ and the original chance constraint (\ref{form:cco-individual}).
In the latter approach, the uncertainties $\xi \sim \Xi$ are partially characterized using directional deviations.
\begin{defn}[Directional Deviations \citep{chen_robust_2007}]
Given a random variable $\xi \in \mathbf{R}$ with zero mean, the forward deviation is defined as
\begin{equation}
	\delta_+(\xi):= \sup_{\theta > 0} \Bigg\{ \sqrt{\frac{2\ln(\mathbb{E}[\exp(\theta \xi)])}{\theta^2}} \Bigg\}
\end{equation}
and the backward deviation is defined as
\begin{equation}
	\delta_-(\xi):= \sup_{\theta > 0} \Bigg\{ \sqrt{\frac{2\ln(\mathbb{E}[\exp(-\theta \xi)])}{\theta^2}} \Bigg\}.
\end{equation}
\end{defn}

\begin{assumption}[\citep{chen_goal-driven_2009}]
\label{assu:convex_support_set}
Let $\mathcal{W}$ denote the smallest closed convex set containing the support $\Xi$ of $\xi$. We assume that the support set is a second-order conic representable set (e.g. polyhedral and ellipsoidal sets).
\end{assumption}
\begin{assumption}[\citep{chen_goal-driven_2009}]
\label{assu:directional_deviation}
Assume the uncertainties $\{\xi_i\}_{i=1}^d$ are zero mean random variables, with a positive definite covariance matrix $\Sigma$. We define the following index set:
\begin{eqnarray}
\mathcal{J}_+ := \{i: \delta_+(\xi_i) < \infty\},~\mathcal{I}_+ := \{i: \delta_+(\xi_i) = \infty\}, \\
\mathcal{J}_- := \{i: \delta_-(\xi_i) < \infty\},~\mathcal{I}_- := \{i: \delta_-(\xi_i) = \infty\}. 
\end{eqnarray} 
\end{assumption}
For notation simplicity, we define two matrices diagonal $P$ and $Q$ as:
\begin{equation*}
	P:= \diag(\delta_+(\xi_1),\cdots,\delta_+(\xi_d)),~Q:= \diag(\delta_-(\xi_1),\cdots,\delta_-(\xi_d)).
\end{equation*}
Major results developed in \citep{chen_robust_2007,chen_goal-driven_2009} are for the individual \emph{linear} chance constraint (\ref{form:cclp-individual}) with decision variables $x_0 \in \mathbf{R}, x \in \mathbf{R}^n$: 
\begin{equation}
\mathbb{P}_\xi\Big( x_0 + \xi^\intercal x \le 0 \Big) \ge 1- \epsilon
\label{form:cclp-individual}
\end{equation}
Its convex approximation using CVaR (or Markov bound) is
\begin{equation}
t + \frac{1}{\epsilon} \mathbb{E}\big[ [x_0 + \xi^\intercal x - t]_+ \big] \le 0 \label{form:cclp-individual-cvar}
\end{equation}
If we are able to find a function $\pi(x_0,x)$ as an upper bound on $\mathbb{E}\big[ [x_0 + \xi^\intercal x]_+ \big]$, then
\begin{equation}
t+\frac{1}{\epsilon}\pi(x_0-t,x) \le 0 \label{form:cclp-individual-cvar-bound}
\end{equation}
is a safe approximation to (\ref{form:cclp-individual-cvar}).
\begin{thm}{\citep{chen_goal-driven_2009}}
\label{thm:bound_cvar}
Suppose that the primitive uncertainty $\xi$ satisfies Assumption \ref{assu:convex_support_set} and \ref{assu:directional_deviation}. The following functions $\pi^i(x_0, x),i=1,\cdots,5$ are upper bounds of $\mathbb{E}_\xi\big[[x_0 + \xi^\intercal x]_+\big]$:
\begin{eqnarray}
\pi^1(x_0,x) &:=& \big[ x_0 + \max_{\xi \in \mathcal{W}} \xi^\intercal x \big]_+ \\
\pi^2(x_0,x) &:=& x_0 + \big[ -x_0 + \max_{\xi \in \mathcal{W}} (-\xi)^\intercal x \big]_+ \\
\pi^3(x_0,x) &:=& \frac{1}{2} \Big(x_0 + \sqrt{x_0^2 + x^\intercal \Sigma x} \Big)
\end{eqnarray}
\begin{eqnarray}
\pi^4(x_0,y) &:=& \inf_{\mu >0} \Bigg\{ \frac{\mu}{\epsilon} \exp\Bigg(\frac{x_0}{\mu}+\frac{u^\intercal u}{2 \mu^2}\Bigg) \Bigg\}.
\end{eqnarray}
where $u_j = \max\{ x_j \delta_+(\xi_j), - x_j\delta_-(\xi_j)\},j=1,\cdots,n$. This bound is finite if and only if $x_j \le 0,~\forall j \in \mathcal{I}_+$ and $x_j \ge 0,~\forall j \in \mathcal{I}_-$.
\begin{eqnarray}
\pi^5(x_0,y) &:=& x_0+\inf_{\mu >0} \Bigg\{ \frac{\mu}{\epsilon} \exp\Bigg(-\frac{x_0}{\mu}+\frac{v^\intercal v}{2 \mu^2}\Bigg) \Bigg\}.
\end{eqnarray}
where $v_j = \max\{ -x_j \delta_+(\xi_j), x_j\delta_-(\xi_j)\},j=1,\cdots,n$. This bound is finite if and only if $x_j \ge 0,~\forall j \in \mathcal{I}_+$ and $x_j \le 0,~\forall j \in \mathcal{I}_-$.
\end{thm}
\begin{rem}
The epigraphs of $\pi^i(x_0,x),~i=1,\cdots,5$ can be represented as second-order cones. Explicit representations depend on the form of $\mathcal{W}$. More details about the representation of (\ref{form:cclp-individual-cvar-bound}) with different choices of $\pi^i(x_0,y)$ can be found in \citep{chen_goal-driven_2009} and \ref{app-sub:cvar_based_convex_approximation}.
\end{rem}

\subsubsection{Constructing Uncertainty Sets} 
\label{ssub:constructing_uncertainty_sets}
We consider the individual linear chance constraint (\ref{form:cclp-individual}) as in Section \ref{ssub:cvar_based_convex_approximation_of_individual_chance_constraints}. The robust counterpart of (\ref{form:cclp-individual}) is
\begin{equation}
\label{form:cclp-individual-rc}
 x_0 + \xi^\intercal x \le 0,~\forall \xi \in \mathcal{U}_\epsilon
\end{equation}
\begin{assumption}
\label{assu:rlp_bounded_support}
$\{\xi_i\}_{i=1}^{d}$ are independent of each other with zero mean and take values on $[-1,1]^d$, i.e. $\mathbb{E}[\xi_i] = 0$ and $\xi_i \in [-1,1]$ for $i=1,2,\cdots,d$.
\end{assumption}
Clearly, under Assumption \ref{assu:rlp_bounded_support}, a natural choice of uncertainty set is the box $\mathcal{U}^{\text{box}} := \{\xi \in \mathbf{R}^d: - \mathbf{1} \le \xi \le \mathbf{1} \}$. Then $\mathcal{F}_{\mathcal{U}}^\text{box} := \{x \in \mathbf{R}^n: f(x,\xi) \le 0,~\forall \xi \in \mathcal{U}^{\text{box}}\}$ is a safe approximation to $\mathcal{F}_\epsilon$, i.e. $\mathcal{F}_{\mathcal{U}}^\text{box} \subseteq \mathcal{F}_\epsilon$. However, using $\mathcal{U}^\text{box}$ leads to $\mathbb{P}(f(x,\xi) \ge 0) = 0 \ll \epsilon$, which causes conservativeness or even infeasibility in many cases. The following choices of uncertainty sets are less conservative.
\begin{thm}[\citep{ben-tal_robust_2009,ben-tal_robust_1999,bertsimas_price_2004}]
\label{thm:rlo_safe_approx}
(\ref{form:cclp-individual-rc}) is a safe approximation to (\ref{form:cclp-individual}) if $\mathcal{U}_\epsilon$ is one of the following:
\begin{subequations}
\begin{align}
& \mathcal{U}_\epsilon^\text{ball}:= \Big\{\xi \in \mathbf{R}^d: \|\xi\|_2 \le \sqrt{2 \ln(1/\epsilon)}\Big\} \\
& \mathcal{U}_\epsilon^\text{ball-box}:= \Big\{\xi \in \mathbf{R}^d: \|\xi\|_\infty \le 1, \|\xi\|_2 \le \sqrt{2 \ln(1/\epsilon)}\Big\} \\
& \mathcal{U}_\epsilon^\text{budget}:= \Big\{\xi \in \mathbf{R}^d: \|\xi\|_1 \le \sqrt{2d \ln(1/\epsilon) }\Big\}
\end{align}
\end{subequations}
And the resulting robust counterparts (RC)s are second-order cone representable (see Chapter 2 of \citep{ben-tal_robust_2009} and \ref{app-sec:representations_of_robust_linear_programs}). 
\end{thm}
It turns out that constructing uncertainty set $\mathcal{U}_\epsilon$ is closely related with the convex approximation framework in Section \ref{ssub:convex_approximation}-\ref{ssub:cvar_based_convex_approximation_of_individual_chance_constraints}.
\begin{thm}[\citep{chen_cvar_2010}]
\label{thm:convex-approximation-uncertainty-set}
Suppose that $\pi(x_0, x)$ is a convex, closed and positively homogeneous, and is an upper bound to $\mathbb{E}_\xi\big[[x_0 + \xi^\intercal x]_+\big]$ with $\pi(x_0,0) = x_0^+$. Then under Assumptions \ref{assu:convex_support_set} and \ref{assu:directional_deviation} and given $\epsilon \in (0,1)$, it holds that for all $(x_0, x)$ such that $\pi(x_0,x) < \infty$, we have 
\begin{equation}
\label{form:cclp-individual-cvar-bound-uncertaintyset}
	\inf_t\Big( t + \frac{1}{\epsilon} \pi(x_0-t, x) \Big) = x_0 + \max_{z \in \mathcal{U}_\epsilon} x^\intercal z
\end{equation}
for some convex uncertainty set $\mathcal{U}_\epsilon$.
\end{thm}
Given an upper bound $\pi(x_0,x)$ on $\mathbb{E}\big[ [x_0 + \xi^\intercal x]_+\big]$ with required properties, the safe approximation (\ref{form:cclp-individual-cvar-bound}) can be represented in the form of $x_0 + \max_{\xi \in \mathcal{U}_\epsilon} \xi^\intercal x$ for some $\mathcal{U}_\epsilon$. Theorem \ref{thm:convex-approximation-uncertainty-set} only proves the existence of a corresponding uncertainty set $\mathcal{U}_\epsilon$. For the $\pi^i(x_0,x)$ functions given in Theorem \ref{thm:bound_cvar}, their corresponding uncertainty sets can be explicitly calculated.
\begin{prop}[\citep{chen_cvar_2010}]
\label{prop:explicit-uncertainty-sets}
For the functions $\pi^i(x_0, y), i=1,2,\cdots,5$ in Theorem \ref{thm:bound_cvar}, their corresponding uncertainty sets are $\mathcal{U}^1_\epsilon \sim \mathcal{U}_\epsilon^5$ below.
\begin{eqnarray}
\mathcal{U}^1_\epsilon &:=& \mathcal{W},\\
\mathcal{U}^2_\epsilon &:=& \bigg\{\xi \in \mathbf{R}^d: \xi= (1-\frac{1}{\epsilon})\zeta,~\text{for some}~\zeta \in \mathcal{W} \bigg\}, \\
\mathcal{U}^3_\epsilon &:=& \bigg\{\xi \in \mathbf{R}^d: \|\Sigma^{-\frac{1}{2}} \xi\|_2 \le \sqrt{\frac{1- \epsilon}{\epsilon}} \bigg \} \\
\mathcal{U}^4_\epsilon &:=& \bigg\{\xi \in \mathbf{R}^d: \exists s,t \in \mathbf{R}^d, \xi=s-t, \nonumber \\
&& \quad \|P^{-1}s + Q^{-1}t\|_2 \le \sqrt{-2\ln(\epsilon)} \bigg \}, \\
\mathcal{U}^5_\epsilon &:=& \bigg\{\xi \in \mathbf{R}^d: \exists s,t \in \mathbf{R}^d,\xi=s-t, \nonumber \\
&&\quad \|P^{-1}s + Q^{-1}t\|_2 \le \frac{1- \epsilon}{\epsilon} \sqrt{2\ln(\frac{1}{1-\epsilon})} \bigg \}.
\end{eqnarray}
where matrices $\Sigma$,$P$ and $Q$ are defined in Assumptions \ref{assu:convex_support_set} and \ref{assu:directional_deviation}.
\end{prop}
Theorem \ref{thm:convex-approximation-uncertainty-set} and Proposition \ref{prop:explicit-uncertainty-sets} demonstrate that the two seemingly different approaches to constructing safe approximations in Section \ref{sub:safe_approximation} are equivalent in many circumstances.

\subsection{Safe Approximation of Joint Chance Constraints} 
\label{sub:safe_approximation_of_joint_chance_constraints}
Although RO has been successful in approximating individual chance constraints, it is rather unsatisfactory in approximating joint chance constraints \citep{chen_cvar_2010}. We restate the joint chance constraint (\ref{form:ccp_joint_cc}) below 
\begin{equation}
	\mathbb{P}_\xi\Big(f(x,\xi) \le 0 \Big) \ge 1 - \epsilon.
\end{equation}
Most RO-based approaches convert a joint chance constraint to several individual chance constraints, then apply the techniques in Section \ref{sub:safe_approximation_of_individual_chance_constraints} on each individual chance constraint. Results along this line are summarized in Section \ref{ssub:converting_a_joint_chance_constraint_to_individual}. Very few approaches directly deal with joint chance constraints, these approaches are mentioned in Section \ref{ssub:other_approaches}.

\subsubsection{Conversion Between Joint Chance Constraints and Individual Chance Constraints} 
\label{ssub:converting_a_joint_chance_constraint_to_individual}
Section \ref{sub:joint_and_individual_chance_constraints} presents two common approaches to converting a joint chance constraint to individual chance constraints. 

First, according to the Boole's inequality or Bonferroni inequality, if $\sum_{i=1}^m \epsilon_i \le \epsilon$, then
the set of $m$ individual chance constraints
\begin{equation}
\label{eqn:individual-cc}
	\mathbb{P}\Big( f_i(x,\xi) \le 0\Big) \le 1 - \epsilon_i,~i=1,\cdots,m
\end{equation}
 is a safe approximation to the joint chance constraint $\mathbb{P}( f(x,\xi) \le 0) \le 1 - \epsilon$. 
The main issue of this approach is the choice of $\{\epsilon_i\}_{i=1}^m$. The problem becomes intractable if taking $\{\epsilon_i\}_{i=1}^m$ as decision variables \citep{nemirovski_convex_2006,chen_cvar_2010}. It remains unclear about how to find the optimal choices of $\{\epsilon_i\}_{i=1}^m$ \footnote{Most people simply choose $\epsilon_i = \epsilon/m$ \citep{nemirovski_convex_2006,chen_robust_2007}, which could be quite conservative if $m$ is a large number.}. Obviously, this approach could be quite conservative in the following two cases: (i) the individual constraints $f_i(x,\xi),~i=1,2,\cdots,m$ are correlated; and (ii) the choices of $\{\epsilon_i\}_{i=1}^m$ are suboptimal. \citep{chen_cvar_2010} provides some deeper observations on the limitation of this approach: the Bonferroni's inequality could still lead to conservativeness even when (i) the individual chance constraints (\ref{eqn:individual-cc}) are independent; and (ii) the optimal choices of $\{\epsilon_i\}_{i=1}^m$ are found. In other words, (\ref{eqn:individual-cc}) is only a safe approximation at best, it may not be equivalent to (\ref{form:ccp_joint_cc}) even with optimal $\{\epsilon_i\}_{i=1}^m$.

The second approach is to define the pointwise maximum of functions $\{f_i(x,\xi)\}_{i=1}^m$ over $x$ and $\xi$, i.e. 
\begin{equation*}
\overline{f}(x,\xi):=\max\Big\{f_1(x,\xi),\cdots,f_m(x,\xi) \Big\}.
\end{equation*}
then the joint chance constraint $\mathbb{P}(f(x,\xi)\le0) \ge 1 - \epsilon$ is equivalent to the individual chance constraint $\mathbb{P}_\xi\big( \overline{f}(x,\xi) \le 0 \big) \ge 1 -\epsilon$. 
The advantage of this approach is that it does not require parameter tuning or induce additional conservativeness. In some cases, e.g. scenario approximation of CVaR in Remark \ref{rem:data_driven_cvar}, this could lead to formulations that are easy to solve \citep{geng_data-driven_2019}. However, in most cases, the structure of $\overline{f}(x,\xi)$ is too complicated to apply the techniques in Section \ref{sub:safe_approximation_of_individual_chance_constraints}.

\subsubsection{Other Approaches} 
\label{ssub:other_approaches}
There might be only three RO-related approaches that directly deal with joint chance constraints. The first approach is robust conic optimization (see Chapter 5-11 of \citep{ben-tal_robust_2009}). The inner constraint $f(x,\xi) \le 0$ is written as a conic inequality, then tractable safe approximations of the robust conic inequality are derived and solved. This approach can model a majority of optimization problems under uncertainties. However, the main limitation is that the resulting robust counterparts are not tractable in many circumstances.

The second approach \citep{chen_cvar_2010} generalizes the CVaR-based convex approximation in Theorem \ref{thm:bound_cvar} and Proposition \ref{prop:explicit-uncertainty-sets}. It proposes a safe approximation to the joint chance constraint (\ref{form:ccp_joint_cc}), and the safe approximation is second-order cone representable. The performance of this approach depends on the choice of a few tuning parameters. Although it is difficult to find the optimal setting, \citep{chen_cvar_2010} designed an algorithm that is guaranteed to improve the choice of parameters. \citep{chen_cvar_2010} also shows that it is possible to combine all the $\pi^i(x_0,x)$ functions in Theorem \ref{thm:bound_cvar} together to reduce conservativeness.
 
The third approach directly dealing with joint chance constraints is the \emph{data-driven robust optimization} proposed in \citep{bertsimas_data-driven_2018}. It shows that by running different hypothesis tests on datasets, it is possible to construct different uncertainty sets that lead to safe approximations of the joint chance constraint (\ref{form:ccp_joint_cc}) with high probability. It is worth noting that the theoretical results in \citep{bertsimas_data-driven_2018} holds for non-convex functions $f(x,\xi)$, albeit the resulting (RC) is very likely to be computationally intractable.

\section{ConvertChanceConstraint (CCC): A Matlab Toolbox} 
\label{sec:convertchanceconstraint_a_matlab_toolbox}
Most existing optimization solvers cannot directly solve (CCO). All reviewed methods in Section \ref{sec:scenario_approach}-\ref{sec:robust_optimization} translate (CCO) to forms that can be recognized and solved by optimization solvers, e.g. SAA converts (CCO) to a mixed integer program (MIP), which can be solved by Gurobi. When solving a chance-constrained program, a typical approach is to write the converted formulation (e.g. the MIP of SAA) in the compact format that a solver recognizes then rely on the solver to get optimal solutions. This approach is unnecessarily repetitive as it needs to be repeated by different researchers on different problems. In addition, different solvers often take various input formats, thus this typical approach is limited to one specific solver. To overcome these issues, an interface or toolbox that automatically converts (CCO) to suitable forms for a variety of solvers is needed.

The remaining part of this subsection introduces the open-source Matlab toolbox \emph{ConvertChanceConstraint} (CCC), which is developed to automate the process of converting chance constraints. CCC is written in Matlab, one of the most popular tools in engineering and many other fields. In consideration of flexibility in modeling and compatibility with existing solvers, CCC is built on YALMIP \citep{lofberg_yalmip_2004}, a modeling language for optimization in Matlab. CCC is open-source on Github \footnote{https://github.com/xb00dx/ConvertChanceConstraint-ccc}, other researchers and engineers could freely use, modify and improve it.

\begin{figure}[tb]
	\centering
	\includegraphics[width=\linewidth]{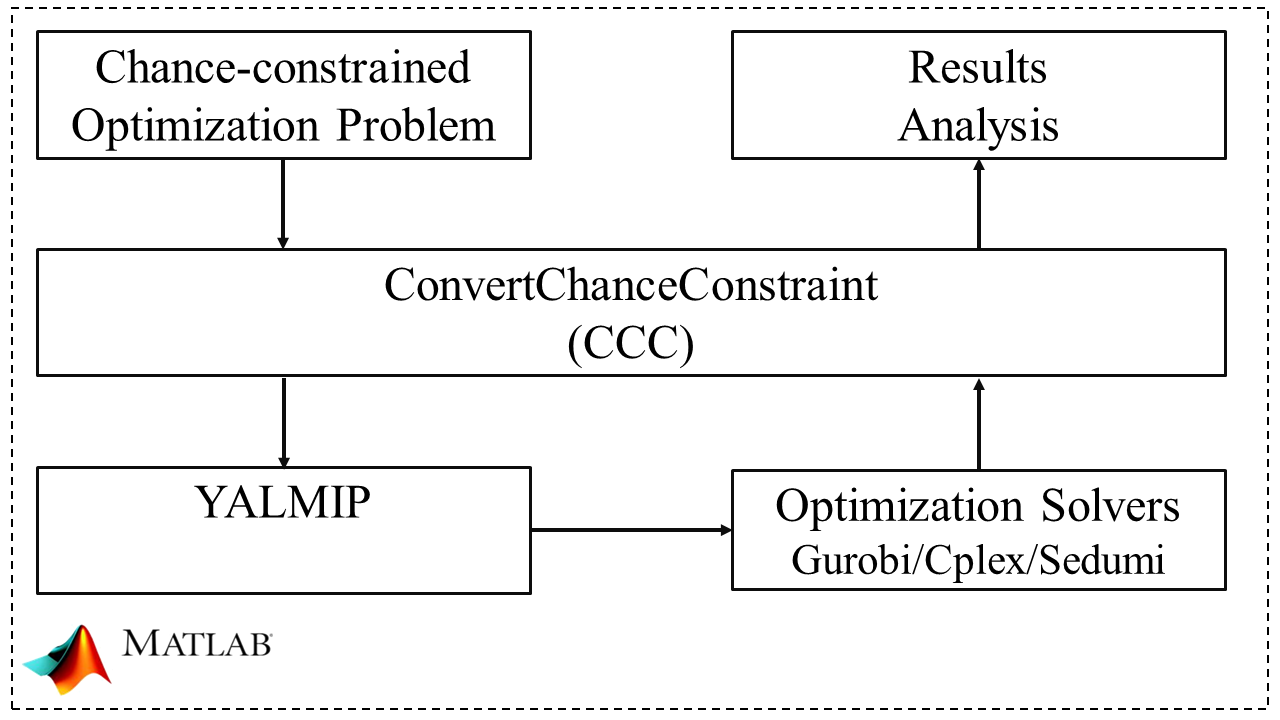}
	\caption{Solving and Analyzing a Chance-Constrained Program via CCC}
	\label{fig:flowchart_ccc}
\end{figure}

Figure \ref{fig:flowchart_ccc} illustrates the logic flow when using CCC to solve and analyze a chance-constrained program. The problem is first formulated in the language of Matlab and YALMIP, then the chance constraint is modeled using the \emph{prob()} function defined in CCC. After receiving the problem formulation and specified method to use (e.g. scenario approach), CCC translates the chance constraint to the formulation that YALMIP could understand. Then YALMIP interfaces with various solvers and further translates the problem for a specific solver. After optimization solver returns the optimal solution, CCC provides a few functions for result analysis, e.g. checking out-of-sample violation probability, calculating the posterior guarantees of the scenario approach.

\begin{figure}[tb]
	\centering
	\includegraphics[width=\linewidth]{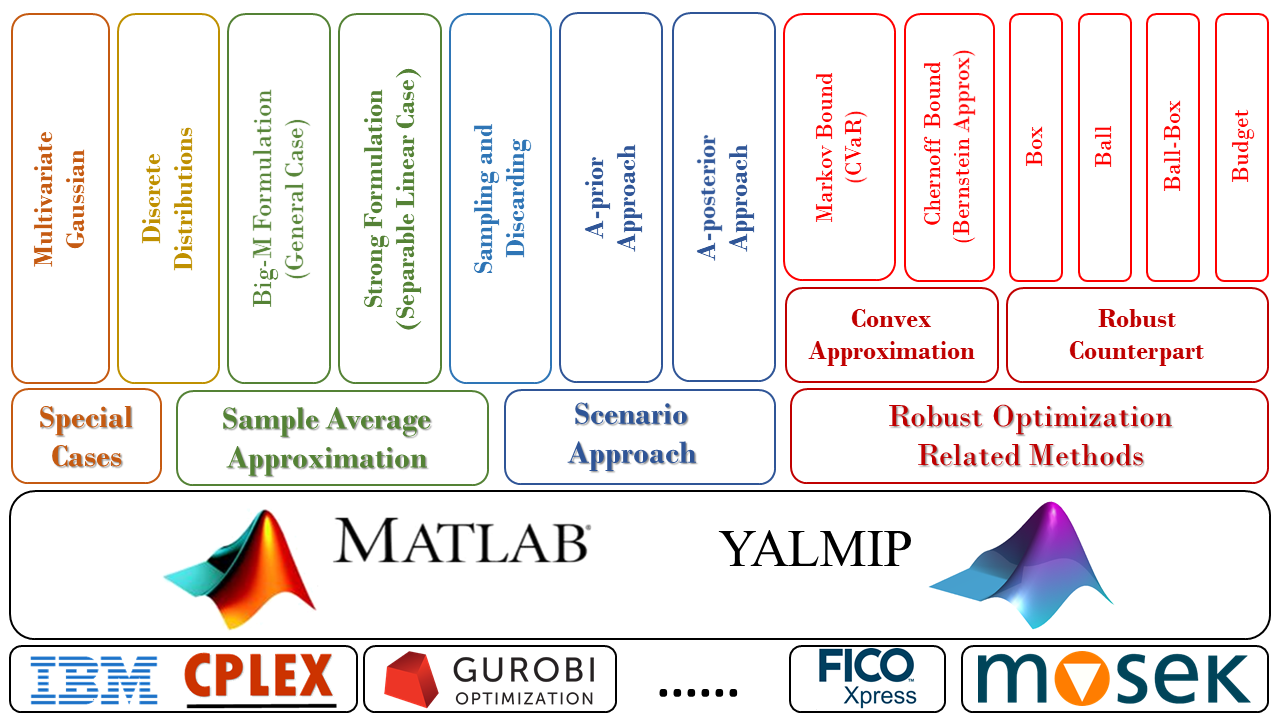}
	\caption{Structure and Main Functions of ConvertChanceConstraint}
	\label{fig:ccc}
\end{figure}

Figure \ref{fig:ccc} presents the structure and main functions of CCC. Three major methods to solve (CCO) are implemented: scenario approach, sample average approximation and robust optimization related methods. The implementation of RO-related methods is based on the robust optimization module \citep{lofberg_automatic_2012} of YALMIP. As illustrated in Figure \ref{fig:flowchart_ccc} and \ref{fig:ccc}, CCC is interfaced via YALMIP with most existing optimization solvers, e.g. Cplex \citep{cplex_v12._2009}, Gurobi \citep{gurobi_optimization_gurobi_2016}, Mosek \citep{mosek_mosek_2015} and Sedumi \citep{sturm_using_1999}.

\section{Concluding Remarks} 
\label{sec:concluding_remarks}
This paper presents a comprehensive review on the fundamental properties, key theoretical results, and three classes of algorithms for chance-constrained optimization. An open-source MATLAB toolbox ConvertChanceConstraint is developed to automate the process of translating chance constraints to compatible forms for mainstream optimization solvers. 

Many interesting directions are open for future research. More thorough and detailed comparisons of solutions to (CCO) on various problems with realistic datasets is needed. 
In terms of theoretical investigation, an analytical comparison of existing solutions to chance-constrained optimization is necessary to substantiate the fundamental insights obtained from numerical simulations.

\appendix

\section{Representations of Robust Linear Programs} 
\label{app-sec:representations_of_robust_linear_programs}


\section{Representations of Convex Approximation} 
\subsection{Sample Approximation of CVaR} 
\label{app-sub:sample_approximation_of_cvar}


\subsection{Representations of $\pi^i(x_0, x)$ Functions in Section \ref{ssub:cvar_based_convex_approximation_of_individual_chance_constraints}} 
\label{app-sub:cvar_based_convex_approximation}


\section*{References}
\bibliographystyle{elsarticle-harv}
\bibliography{references}

\end{document}